# Application of an accurate remainder term in the calculation of Residue Class Distributions

COLIN MYERSCOUGH

ABSTRACT. This paper uses concepts introduced by Fiorilli and Martin to develop a more accurate remainder term in calculations of the distribution of primes in residue classes. Using it in Rubinstein and Sarnak's method estimates the limiting logarithmic frequency of $\pi(x) > \text{Li}(x)$ as $2.6300 \times 10^{-7}$ with only 5 zeta zeros used explicitly, and a reliable value between $2.62996732 \times 10^{-7}$ and $2.62996733 \times 10^{-7}$ is obtained with fewer than 100 zeros. Accurate results for "prime number races" can usually be obtained with explicit use of only the first zero of each $L$-function involved, thus bringing out their dependence on those zeros. For extreme deviations, the method of steepest descent can be applied, and the remainder approximated by an explicit formula together with a rapidly convergent series. This gives, for example, $8.649 \times 10^{-12477}$ for the logarithmic frequency of $\pi(x) > \text{Li}(x) + 5\text{Li}(\sqrt{x})$. The Monach - Lamzouri model of the extreme distribution is developed to give close agreement with these results. The remainder can also be calculated explicitly as an asymptotic series. This allows good modelling of non-extreme races and the calculation of distributions by direct convolution. The three different methods agree to within 0.001% where they overlap in application. The distribution calculated from values of $\pi(x)$ for $x \leq 10^{20}$ shows similar behaviour to the limiting distribution, but is somewhat further from normal.

## 1. INTRODUCTION

As set out elsewhere, notably in Rubinstein and Sarnak 1994, Feuerverger and Martin 2000, and Fiorilli and Martin 2012, and summarised in <u>Section 3</u> below, assuming the Riemann hypothesis (generalised as necessary), the limiting relative logarithmic distribution of primes between residue classes $a$ and $b$ mod $q$, or of all primes, has density $P_0(v)$, where $P_u(v)$ is the density of

$$2 \sum_{\chi} \alpha_\chi \sum_{\gamma > u} \frac{X_\gamma}{\sqrt{1/4 + \gamma^2}} \quad \text{where} \quad 0 \leq \alpha_\chi = \frac{1}{2}|\chi(a) - \chi(b)| \leq 1 \qquad (1.1)$$

The sum in (1.1) runs over the non-trivial Dirichlet characters $\chi$ determining the distribution, and positive imaginary parts $\gamma$ of zeros of the corresponding $L$-functions, which zeros are assumed to be linearly independent in the rational field. If the limiting logarithmic distribution of the difference between squares and non-squares mod $q$, or of normalised fluctuation of prime count, is being considered, there is just one $L$– or zeta function involved, and $\alpha_\chi = 1$. The $X_\gamma$ are independent random variables, each with the density:

$$f(x) = \frac{1}{\pi\sqrt{1 - x^2}} \quad (|x| < 1), \qquad 0 \ (|x| > 1) \qquad (1.2)$$

All methods of calculating $P_0(v)$, and the probability of normalised deviation exceeding $v$:

$$E(v) = \int_v^\infty P_0(v')dv' \tag{1.3}$$

depend on combining an approximation for $P_u(v)$ with the exact effect for each $\chi$ of the $N_\chi(u)$ separate $X_\gamma$ with $\gamma < u$. With the approximations previously used, accurate results require values of many thousands for $N_\chi$. This paper develops a much better approximation for $P_u(v)$, of near normal (Gaussian) form. Using this, the computational requirements of Rubinstein and Sarnak's method are reduced by a factor of $10^5$ or more, to give accurate results with $N_\chi \leq 10$, and to the limit of computational accuracy with $N_\chi = 100$. Results for large $v$ may be calculated accurately for the first time, by other methods presented here. The results in Table 1 were worked out on an ordinary desktop micro using Microsoft Excel.

**TABLE 1.** Limiting logarithmic frequency of $\pi(x) > \text{Li}(x) + m\text{Li}(\sqrt{x})$

| m | Value | $N_\chi$ | Method |
|---|---|---|---|
| 0 | $2.6300 \times 10^{-7}$ | 5 | Rubinstein - Sarnak |
| 0 | range $2.62996732 \times 10^{-7}$ to $2.62996733 \times 10^{-7}$ | 100 | Rubinstein - Sarnak |
| 1 | $1.603 \times 10^{-95}$ | 250 | Steepest descent/convolution |
| 2 | $6.239 \times 10^{-489}$ | 1500 | Steepest descent |
| 3 | $7.901 \times 10^{-1724}$ | 3000 | Steepest descent |
| 4 | $5.621 \times 10^{-4954}$ | 10000 | Steepest descent |
| 5 | $8.649 \times 10^{-12477}$ | 16000 | Steepest descent |

## 2. STATEMENT OF KEY RESULTS

**Preliminaries and notation.** The process of calculating $P_0(v)$ may be expressed as a convolution

$$P_0(v) = \int_{-S(u)}^{S(u)} F_u(v') P_u(v - v') dv' \tag{2.1}$$

where $F_u(v)$ is the probability density of the sum $2 \sum_\chi \alpha_\chi \sum_{0 < \gamma < u} \frac{X_\gamma}{\sqrt{1/4 + \gamma^2}}$, which can be obtained by repeated convolution of densities $\frac{2\alpha_\chi}{\sqrt{1/4+\gamma^2}} f\left(\frac{v\sqrt{1/4+\gamma^2}}{2\alpha_\chi}\right)$, and

$$S(u) = \sum_\chi \alpha_\chi S_\chi(u) \quad \text{where} \quad S_\chi(u) = 2 \sum_{0 < \gamma < u} \frac{1}{\sqrt{1/4 + \gamma^2}} \tag{2.2}$$

is the 'span' of $F_u(v)$; $F_u(v) = 0$ for $|v| > S(u)$. The number of initial terms is

$$N(u) = \sum_\chi N_\chi(u) \quad \text{where for each } \chi \quad N_\chi(u) = \sum_{0 < \gamma < u} 1 \tag{2.3}$$

The convolutions are generally done using Fourier Transforms (Rubinstein-Sarnak) or 2 sided Laplace Transforms (method of steepest descent, as applied by Monach 1980 and Lamzouri 2012 to establish the behaviour of $P_0(v)$ for large $v$, though not to computations.) Using throughout this paper the representation of such transforms for $P_0(v)$ or any other function:

$$\hat{P}_0(\omega) = \int_{-\infty}^{\infty} P_0(v) e^{-i\omega v} dv, \qquad \bar{P}_0(s) = \int_{-\infty}^{\infty} P_0(v) e^{sv} dv, \qquad (2.4)$$

then

$$\hat{P}_0(\omega) = \hat{F}_u(\omega)\, \hat{P}_u(\omega); \qquad \bar{P}_0(s) = \bar{F}_u(s)\bar{P}_u(s). \qquad (2.5)$$

$P_u(v)$ and other functions considered here are symmetric, therefore their Fourier Transforms are real. Since $\hat{f}(\omega) = \frac{1}{\pi}\int_{-1}^{1} \frac{e^{-i\omega v} dv}{\sqrt{1-v^2}} = \frac{1}{\pi}\int_0^{\pi} \cos(\omega \cos\theta)\, d\theta = J_0(\omega)$, a Bessel function, and similarly $\bar{f}(s) = I_0(s)$, a Bessel function of imaginary argument:

$$\hat{F}_u(\omega) = \prod_{\chi} \prod_{0<\gamma<u} J_0\left(\frac{2\alpha_\chi \omega}{\sqrt{1/4+\gamma^2}}\right) \qquad \bar{F}_u(s) = \prod_{\chi} \prod_{0<\gamma<u} I_0\left(\frac{2\alpha_\chi s}{\sqrt{1/4+\gamma^2}}\right) \qquad (2.6)$$

$$\hat{P}_u(\omega) = \prod_{\chi} \prod_{\gamma>u} J_0\left(\frac{2\alpha_\chi \omega}{\sqrt{1/4+\gamma^2}}\right) \qquad \bar{P}_u(s) = \prod_{\chi} \prod_{\gamma>u} I_0\left(\frac{2\alpha_\chi s}{\sqrt{1/4+\gamma^2}}\right). \qquad (2.7)$$

Also define

$$L_u(s) = \log \bar{P}_u(s) = \sum_{\chi} \sum_{\gamma>u} \log I_0\left(\frac{2\alpha_\chi s}{\sqrt{1/4+\gamma^2}}\right) \qquad (2.8)$$

Thus $\hat{P}_u(\omega)$, $\bar{P}_u(s)$ and $L_u(s)$ are respectively the characteristic function, moment generating function, and cumulant generating function of $P_u(v)$. All these functions and distributions will sometimes be discussed for individual Dirichlet characters for which they have no direct meaning in terms of a logarithmic distribution. Also define

$$b_k(\chi, u) = \sum_{\gamma>u} \frac{1}{(1/4+\gamma^2)^k}; \qquad r_k(\chi,u) = \frac{b_k(\chi,u)}{(b_1(\chi,u))^k} \qquad (2.9)$$

$$B_k(u) = \sum_{\chi} \alpha_\chi^{2k}\, b_k(\chi,u); \qquad R_k(u) = \frac{B_k(u)}{\left(B_1(u)\right)^k} \qquad (2.10)$$

so that the standard deviation $\sigma_u$ of $P_u(v)$ is $\sqrt{2B_1(u)}$. The dependence of these quantities on $\chi$ and $u$ is omitted below when no ambiguity is thereby generated, to simplify notation.

**Behaviour of $P_u(v)$ for small $v$.** This is discussed initially in <u>Section 4</u>. As $u$ increases, $P_u(v)$ appears to approach normal form over a wider and wider range of $Y = v/\sigma_u$. This approach is quantified by expressing in terms of the $R_k(u)$ the deviation of the moments of the distribution from those pertaining to a normal distribution, standard deviation $\sigma_u$. The behaviour of $R_k(u)$ is given through the following results, which are more accurate than those previously derived.

**Theorem 1.** *For any Dirichlet character $\chi$, and integer $k \geq 1$,*

$$b_k = \frac{y + 1/(2k-1)}{2(2k-1)\pi u^{2k-1}}\left[1 + O\left(\frac{k}{u}\right)\right] \qquad (2.11)$$

$$r_2 = \frac{2\pi(y+1/3)}{3u(y+1)^2}\left[1 + O\left(\frac{1}{u}\right)\right] \tag{2.12}$$

$$\frac{r_k}{r_2^{k-1}} = \frac{y+1/(2k-1)}{(2k-1)(y+1)}\left[\frac{3(y+1)}{(y+1/3)}\right]^{k-1}\left[1 + O\left(\frac{k}{u}\right)\right] \tag{2.13}$$

$$S_\chi - \frac{\log u}{\pi}\left[\frac{\log u}{2} + \log\left(\frac{q^*}{2\pi}\right)\right] = \Delta_\chi + O\left(\frac{\log u}{u}\right) \tag{2.14}$$

where $y = \log\left(\frac{q^*u}{2\pi}\right)$, $q^* \leq q$ is the conductor ($q^* = q$ for all $\chi$ when $q$ is an odd prime), and the constant $\Delta_\chi$ depends largely on the values of the smallest $\gamma$.

**Theorem 2.** *For large $u$, $R_k(u)/(R_2(u))^{k-1}$ is reduced below the average value of $r_k(\chi,u)/(r_2(\chi,u))^{k-1}$ by a factor of at least $(\phi(q)/2)^{k-2}$, $\phi(q)$ being Euler's totient function.*

$\beta = R_2(0)$ (values for some distributions in Table 2 below) is the main determinant of the difference of $P_0(v)$ from a normal distribution at small $v$. An inter-residue distribution dependent on several characters will be closer to normal than any of the distributions determined by single characters.

**Behaviour of $P_u(v)$ for large $v$.** However, consideration of moments also shows that for any $u$, $P_u(v)$ will differ significantly from normal for large enough $v$. In Section 5, an argument given in Montgomery 1979 and Montgomery and Odzylko 1988 is extended to prove

**Theorem 3.** *For any $u$, when $v \geq S(u)$*

$$\log E(v) < -\frac{(v - S(u))^2}{2\sigma_u^2} \tag{2.15}$$

*and consequently $P_0(v)$ and $E(v)$ tend to zero as $v \to \infty$ faster than any normal distribution.*

**Difference of $\hat{P}_u(\omega)$ from transform of normal distribution.** In Section 6 this is accurately quantified through

**Theorem 4.** *Defining $\tau = \omega\sigma_u$, for $|\tau| < T$ and $K \geq 1$*

$$\hat{P}_u(\omega) = \exp\left(-\sum_{k=1}^{K} c_k R_k \tau^{2k}\right)\left[1 + O\left(\frac{c_K R_K \tau^{2K+2}}{(T^2 - \tau^2)}\right)\right] \tag{2.16}$$

Here $c_1 = 1/2$, $c_2 = 1/16$, $c_3 = 1/72$, $c_4 = 11/3072$, $c_5 = 19/19200$, whilst for $k > 5$

$$\frac{1}{k}\left(\frac{2}{j_1^2}\right)^k < c_k < \frac{1}{k}\left(\frac{2}{j_1^2}\right)^k[1 + 1.16(0.19)^k] \tag{2.17}$$

where $j_1 \cong 2.2048$ is the first zero of $J_0(\omega)$. For single and multiple series respectively

$$T = j_1\sqrt{\frac{y+1/3}{3(y+1)r_2}} \quad ; \quad T \approx j_1\sqrt{\frac{(y+1/3)\phi(q)}{6(y+1)R_2}} \tag{2.18}$$

**Application to Rubinstein and Sarnak's method.** The results for $\pi(x) > \text{Li}(x)$ ($m = 0$) in Table 1 for were obtained by applying (2.16) with $K = 7$ in this method. Calculations can be further accelerated by using (2.15) to validate use of 10 to 20 terms, rather than typically 1000 in previous work, in the Poisson summation employed to calculate $P_0(v)$ or $E(v)$ from values of $\hat{P}_0(\omega)$. The only practical limitations on accuracy obtainable are those of data on zeros used, and of computation (about $10^{-16}$ in the work presented here).

TABLE 2. PRIME NUMBER RACE RESULTS

| q | Race | Chars | $\sigma_0 = \sqrt{2B_1(0)}$ | $\beta = R_2(0)$ | $N_\chi(u)$ | K | result |
|---|---|---|---|---|---|---|---|
| 1 | Pi vs. Li | 1 | 0.2149 | 0.0696 | 5 | 7 | $2.6300 \times 10^{-7}$ |
| 4 | 1 leads 3 | 1 | 0.3944 | 0.1517 | 3 | 5 | 0.0040721 |
| 5 | square leads | 1 | 0.3957 | 0.1161 | 3 | 5 | 0.0045774 |
| 5 | 1 leads 2 or 3 | 3 | 0.3598 | 0.0571 | 1 | 5 | 0.0478254 |
| 7 | square leads | 1 | 0.5052 | 0.1871 | 3 | 5 | 0.0217412 |
| 7 | 1 leads 3 or 5 | 5 | 0.8666 | 0.0493 | 1 | 4 | 0.1255461 |
| 7 | 1 leads 6 | 3 | 0.9658 | 0.1316 | 1 | 4 | 0.1547904 |
| 8 | 1 leads 3 | 2 | 0.6253 | 0.1926 | 2 | 5 | 0.0004312 |
| 8 | order 7, 1, 5 | 3 | | | 5 | 3 | 0.0024769 |
| 13 | square leads | 1 | 0.6298 | 0.2741 | 3 | 5 | 0.0556810 |
| 13 | 1 leads 6 or 11 | 11 | 1.6570 | 0.0427 | 1 | 4 | 0.2745797 |
| 13 | 1 leads 5 or 8 | 9 | 1.8185 | 0.1022 | 1 | 3 | 0.2953206 |
| 13 | 1 leads 2 or 7 | 11 | 2.0384 | 0.2020 | 1 | 3 | 0.3211909 |

Table 2 gives some results for "prime number races" (given as logarithmic probability of the less common residue leading and therefore equal to $E(2^l)$ for certain $l \geq 0$; $l = 0$ when $q$ is an odd prime). They agree with other published calculations to the accuracy of those. For prime modulus inter-residue races, precision to 7 decimal places can usually be obtained with $N_\chi = 1$ (instead of typically 10,000 in previous calculations), thus bringing out the dependence of these results on the smallest imaginary parts of zeros, as previously observed by Bays *et al.* 2001.

**Method of steepest descent.** This (in the form originally introduced by Laplace) depends on the variation of $P_0(v)e^{sv}$ near its maximum being close to a normal distribution. Section 7 describes how to use this to calculate small values of $P_0(v)$ for large $v$. Rubinstein and Sarnak's method is not suitable for this purpose. To make such calculations accurately requires a clear understanding of the behaviour of $L_0(s)$ and its derivatives for large $s$. This is given by

**Theorem 5.** *For large s and a single character distribution*

$$L_0'(s) = \frac{(\log s)^2}{2\pi} + \frac{A \log s}{\pi} + (\Delta_\chi + X) + O\left(\frac{\log s}{s}\right);$$

$$L_0(s) = sL_0'(s) - \frac{s}{\pi}\{\log s + A - 1\} + O(\log s)$$

$$L_0''(s) = \frac{(\log s + A)}{\pi s} + O\left(\frac{\log s}{s^2}\right); \quad L_0'''(s) = \frac{-(\log s + A - 1)}{\pi s^2} + O\left(\frac{\log s}{s^3}\right) \quad (2.19)$$

Here

$$A = A_0 + \log\left(\frac{q^*}{\pi}\right) \text{ where } A_0 = 1 + \int_0^1 \frac{\log I_0(x)}{x^2} dx + \int_1^\infty \frac{\log I_0(x) - x}{x^2} dx \cong -0.08933 \quad (2.20)$$

$$A_1 = \int_0^1 \log x \frac{\log I_0(x)}{x^2} dx + \int_1^\infty \log x \frac{(\log I_0(x) - x)}{x^2} dx \cong -2.12634 \quad (2.21)$$

$$X = \frac{1}{\pi}\left[(\log 2 + A_0)\log\left(\frac{q^*}{\pi}\right) - \frac{(\log 2)^2}{2} - 1 + A_0 - A_1\right] \cong 0.0336 + 0.1922 \log q^* \quad (2.22)$$

The proof follows the methods of Monach 1980 and Lamzouri 2012, but is taken to higher accuracy. $A_0$ is the same quantity as referred to by them, though note its correct value.

**Calculations using method of steepest descent.** In <u>Section 8</u> using Theorem 5 is proved the following, accurate to third order in the parameter $\kappa(s)$ which reflects the extent to which the variation of $P_0(v)e^{sv}$ near its maximum deviates from a normal distribution:

**Theorem 6.** *For $V = L_0'(s)$ :*

$$\log P_0(V) = L_0(s) - sV - \tfrac{1}{2}\log(2\pi L_0''(s)) + \kappa^2\left(\tfrac{\lambda}{8} + \tfrac{1}{6}\right) + O((s \log s)^{-2}) \quad (2.23)$$

*where $\kappa^2 \sim \frac{\pi}{s \log s}$ and $\lambda \sim -1$ can be calculated in terms of derivatives of $L_0(s)$.*

<u>Section 9</u> describes how using Theorem 1 $L_u(s)$ and its derivatives, and hence

$$L_0(s) = \sum_\chi \sum_{0 < \gamma < u} \log I_0\left(\frac{2\alpha_\chi s}{\sqrt{1/4 + \gamma^2}}\right) + L_u(s) \quad (2.24)$$

and its derivatives, may be calculated with high accuracy, The result corresponding to Theorem 4 is

**Theorem 7.** *Defining $t = \sqrt{2B_1 s}$, for $|t| <$ T as defined by (2.18) and $K \geq 1$*

$$L_u(s) = \sum_{k=1}^{K} (-1)^{k-1} c_k R_k t^{2k} + O\left(\frac{c_K R_K t^{2K+2}}{(T^2 - t^2)}\right) \quad (2.25)$$

*which expression may be differentiated to obtain derivatives of $L_u(s)$.*

(2.24) and (2.25) allow values of $L_0(s)$ and its derivatives to be calculated, and (2.23) then gives $\log P_0(V)$ for $V = L_0'(s)$. However, the calculations demand high accuracy. For the prime count distribution at $V = 11$, (2.19) shows (and (2.24) as well as calculations verifies) that $s$ is of order 10000, $L_0(s)$ of order 100000, and $\log P_0(V)$ of order -29000. If $P_0(V)$ is to be calculated to a particular accuracy, say 0.01%, values of $L_0(s)$ are required to accuracy better than 0.0001, or 1 part in $10^9$. The main contribution to $L_0'(s)$ is from $S_\chi$ or combinations thereof, and then requires $N_\chi(u) \approx 15000 - 20000$. To limit $N_\chi(u)$ requires calculation to within 0.0001 of values of $L_u(s)$ which are between 10000 and 20000, for values of $t/T$ which are significant fractions of 1. Direct summation of the series $\sum_{k=1}^\infty (-1)^{k-1} c_k R_k t^{2k}$ would need many terms to achieve this. However, summation can be accelerated through

**Theorem 8.** *For q prime, and $|t| < T$, $L_u(s) = \Sigma(t) + \Sigma_E(t)$ where*

$$\Sigma(t) = \frac{T^2}{(y+1)j_1^2}\left[(y-1)\left(\frac{2t}{T}\arctan\left(\frac{t}{T}\right) - \log\left(1 + \left(\frac{t}{T}\right)^2\right)\right) + \frac{2t}{T}\int_0^{t/T} \frac{\arctan w\, dw}{w}\right] \quad (2.26)$$

*and the terms of the series for $\Sigma_E(t)$ eventually decrease more quickly than $(0.19\, t^2/T^2)^k$*

(A similar, but more difficult to express, result is available for general $q$.)

To a first approximation, $\log E(v) = \log P_0(v) - \log s$. A more accurate approximation is derived by methods similar to the derivation of (2.23). Table 4 gives $\log E(v)$ for values of $v$ up to 11, which correspond to the event $\pi(x) > \text{Li}(x) + 5\text{Li}(\sqrt{x})$ and to similar events for other races. The accuracy of the method improves as $v$ increases. For the prime count distribution, 4 significant figures (0.01% proportional accuracy in $E(v)$) is attained for $v < 1.5$.

A more accurate model than previously given of the 'double exponential' behaviour of $P_0(v)$ and $E(v)$ for large $v$, which throughout the ranges considered agrees with calculations to within a factor of less than 10 ($\log P_0(v)$ and $\log E(v)$ accurate to within 1 – 2), is given by

**Theorem 9.** *For prime $q$*

$$\log P_0(v) = -\frac{\alpha_\Sigma}{q}(W-1)e^{W-Y-A_0} + O\left(\sqrt{2\pi v/\alpha_\Sigma}\right) \quad (2.27)$$

*and $\log E(v)$ is similarly expressed to this accuracy, where*

$$W^2 = A^2 + Y^2 + \frac{2\pi v}{\alpha_\Sigma} - \frac{2\pi}{\alpha_\Sigma}\sum_\chi \alpha_\chi \Delta_\chi - 2X - Z$$

$$\alpha_\Sigma = \sum_\chi \alpha_\chi \quad Y = \frac{1}{\alpha_\Sigma}\sum_\chi \alpha_\chi \log \alpha_\chi \quad Z = \frac{1}{\alpha_\Sigma}\sum_\chi \alpha_\chi (\log \alpha_\chi)^2 \quad (2.28)$$

*and X is defined by (2.22).*

A similar result is available for general $q$, but is more difficult to express.

**Explicit expression for $P_u(v)$.** Section 10 develops this, proving

**Theorem 10.** *Defining $Y = v/\sigma_u$, and $T$ by (2.18), for $Y$ rather less than $T$*

$$P_u(v) = \frac{G_u(Y)e^{-Y^2/2}}{\sigma_u\sqrt{2\pi}} \quad (2.29)$$

*where $G_u(Y)$ can be represented by the asymptotic expansion*

$$G_u(Y) \approx \sum_{m=0}^{M} \frac{(-1)^m}{2^m m!}\frac{d^{2m}}{dY^{2m}}\exp\left(\sum_{k=2}^{K}(-1)^{k-1}c_k R_k Y^{2k}\right) \quad (2.30)$$

Thus $G_u(Y)$ describes explicitly the difference of $P_u(v)$ from normal form. For an inter-residue race when $q$ is large (making $\sigma_0$ large and $R_2$ small as explained in Section 3) $G_0(v/\sigma_0)$ represents well this difference for $P_0(v)$ in $0 \leq v \leq 1$. Then the race result for 1 leading a non-square is:

$$E(1) \cong \frac{1}{2} - \frac{1}{\sigma_0 \sqrt{2\pi}} \left(1 - \frac{3\beta}{16} - \frac{1}{6\sigma_0^2}\right). \tag{2.31}$$

which is equivalent to the unnumbered equation on page 40 of Fiorilli and Martin 2012.

The explicit expression for $P_u(v)$ also allows the calculation of distributions by numerical convolution, as described in <u>Section 11</u>. This involves more computation than using transforms, but delivers the whole distribution in one calculation and has particular application to multi-way races with real characters ($q = 8, 12$ or $24$). It can be shown that the result of the square/non-square race mod 24 (residue 1 vs. all others) is less than $2.86 \times 10^{-8}$.

**Comparison of results, and conclusions.** In <u>Section 12</u> it is shown that for the prime count distribution, results from the Rubinstein-Sarnak, steepest descent, and convolution methods agree to within 0.001% where they overlap in application. Since the steepest descent method becomes more accurate with increasing $v$, this supports the accuracy of its results for larger $v$. The comparison with actual distributions calculated from values of $\pi(x)$ for $x \leq 10^{20}$ is discussed; the actual distributions show are somewhat further from normal, showing that limiting behaviour is approached quite slowly. Results in the range $v \leq 1.5$ are shown to be consistent with estimates of first excursions of $\pi(x)$. The link between 'near normal' behaviour for small $v$, and 'double exponential' behaviour for large $v$, is summarised.

## 3. SUMMARY OF BASIC THEORY

The history of the study of 'bias' in prime number races goes back to Chebyshev's insights of 1853, but as shown by Rubinstein and Sarnak the key results arise from considering the logarithmic density of excess of one residue class over another. See the introductory sections of Feuerverger and Martin 2000, and Fiorilli and Martin 2012, and the work leading to Proposition 2.6 of the latter, for an account of this; in this paper, some differences in notation have proved convenient. $P_u(v)$ is the limiting distribution of values $v$ taken by the series

$$D_u(z) = -2 \sum_\chi \alpha_\chi \sum_{\gamma > u} \frac{\sin(\gamma z + \operatorname{arc cot}(2\gamma))}{\sqrt{1/4 + \gamma^2}} \tag{3.1}$$

which when $u = 0$ determines the relative behaviour of residue class counts at large $x = e^z$. Thus for the prime count distribution $D(x) = D_0(z) + O(1/\sqrt[4]{x})$ where

$$D(x) = [\pi_0(x) - \operatorname{Li}(x)] \frac{\log x}{\sqrt{x}} \tag{3.2}$$

The existence of $P_u(v)$ and $E(v)$ is discussed in Section 4 below but is not considered in detail here, as it follows from Theorem 1.1 of Rubinstein and Sarnak 1994.

In an inter-residue race mod $q$ there are $\phi(q) - 1$ characters, excluding the trivial character which is left out of summations, and

$$\sum_{\chi} \alpha_\chi^2 = \frac{\phi(q)}{2}; \quad \sum_{\chi} \alpha_\chi^{2k} \leq \sum_{\chi} \alpha_\chi^{2(k-1)} \quad (k \geq 2) \tag{3.3}$$

## 4. BEHAVIOUR OF $P_u(v)$ FOR SMALL $v$

FIGURE 1. PRIME COUNT DENSITY

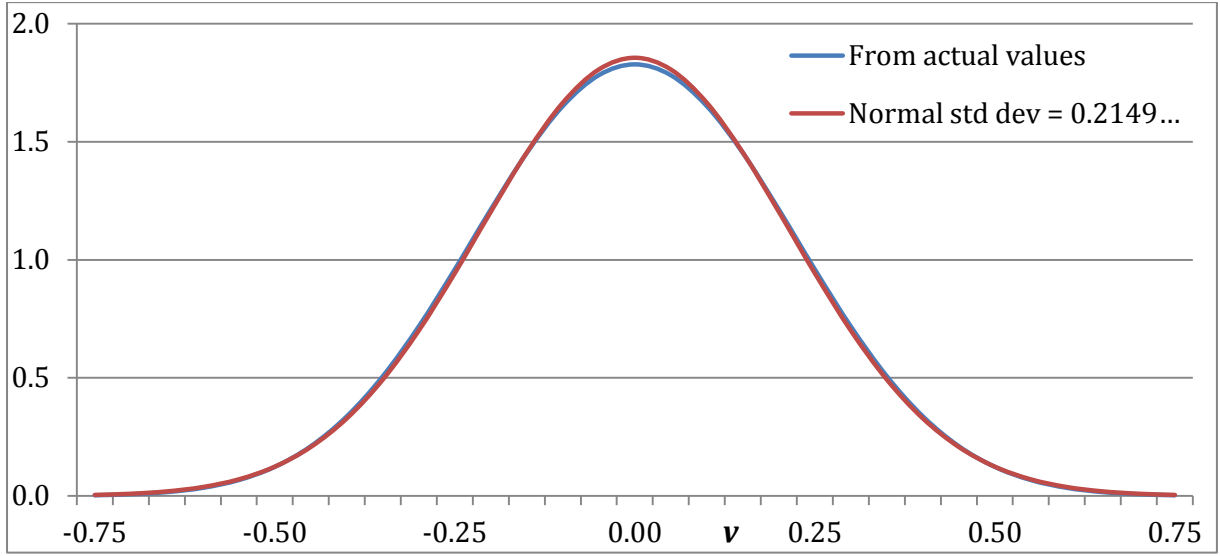

**The near normal distribution.** Figure 1 shows a histogram of $D(x)$ as defined by (3.2), for some 35,000 values of $z = \log x$ spaced at interval 0.001 between 11.513 and 46.051 ($x = 10^5$ to $10^{20}$: values obtained by interpolation from the tables of Kulsha 2011.) Also shown is a normal distribution, standard deviation $\sigma_0 \cong 0.2149$. In the range $|v| < 0.8$ for which the histogram of $D$ is meaningful, its density is slightly lower than normal in the centre, slightly higher at 1 - 2 standard deviations, and lower further out. (This is the type of behaviour shown by (2.30), as is discussed in more detail in Section 10). The distribution differs from normal because the densities combined extend over different ranges $\frac{\sqrt{1/4+\gamma^2}}{2\alpha_\chi}$, as well as being of the far from normal "cup shaped" nature of (1.2). One would expect the difference from normal to be larger if the differences amongst these ranges are greater.

**Results on moments.** Noting that for any zeta or $L$-function, the proportionate variation of successive $\gamma$ decreases with increasing magnitude, it can be expected that the density $P_u(v)$ would become closer to normal form as $u$ increases. This may be quantified by considering its even moments

$$M_{2k} = \int_{-\infty}^{\infty} v^{2k} P_u(v) dv = \lim_{Z \to \infty} \frac{2^{2k}}{Z} \int_{Z'}^{Z'+Z} \left[ \sum_{\chi} \alpha_\chi \sum_{\gamma>0} \frac{\sin(\gamma z + \operatorname{arc cot}(2\gamma))}{\sqrt{1/4+\gamma^2}} \right]^{2k} dz \tag{4.1}$$

Thus

$$\sum_{k=0}^{\infty} \frac{M_{2k} s^{2k}}{(2k)!} = \int_{-\infty}^{\infty} P_u(v) \sum_{k=0}^{\infty} \frac{v^{2k} s^{2k}}{(2k)!} dv = \bar{P}_u(s) = \prod_{\chi} \prod_{\gamma>u} I_0 \left( \frac{2\alpha_\chi s}{\sqrt{1/4+\gamma^2}} \right) \tag{4.2}$$

Each coefficient with $k > 1$ in the series $I_0(2w) = \sum_{k=0}^{\infty} \frac{w^{2k}}{(k!)^2}$ is less than the corresponding coefficient in the series $\exp(w^2) = \sum_{k=0}^{\infty} \frac{w^{2k}}{k!}$. Therefore for $k > 1$ $\frac{M_{2k}}{(2k)!}$, being the sum of products of such terms, is less than the $s^{2k}$ coefficient in the product

$$\prod_{\chi} \prod_{\gamma > u} \exp\left(\frac{\alpha_\chi^2 s^2}{1/4 + \gamma^2}\right) = \exp(B_1(u)s^2) \tag{4.3}$$

which coefficient is $\frac{B_1^k}{k!} = \frac{M_{2k}^{(0)}}{(2k)!}$, where $M_{2k}^{(0)} = \frac{(2k)!B_1^k}{k!}$ is the $2k$th moment of a normal distribution variance $\sigma_u^2$. So $M_{2k} < M_{2k}^{(0)}$. This also confirms that the series (4.2) and product (4.3) converge: $\bar{P}_u(s)$ and $L_u(s)$ exist for all $s$ (and $\hat{P}_u(\omega)$ for all $\omega$). This is sufficient to confirm the existence of $P_u(v)$.

Since the product $v^{2k} \exp(-v^2/2\sigma_u^2)$ has its maximum at $v = \sigma_u\sqrt{2k}$, the ratios $M_{2k}/M_{2k}^{(0)}$ should give an indication of the extent of difference of $P_u(v)$ from normal form near these values of $v$. In this Section and the next will be proved Lemma 1 below, which suggests that as $u$ increases $P_u(\Upsilon\sigma_u)$ approaches its normal value for any fixed $\Upsilon$, the approach requiring larger $u$ for larger $\Upsilon$; but as $\Upsilon$ increases $P_u(\Upsilon\sigma_u)$ will differ increasingly from its normal value for any fixed $u$, for larger $u$ the difference requiring larger $\Upsilon$ to become established. In later Sections of this paper these differences are quantified and applied to calculations.

**Lemma 1.** *For any fixed $u$, $M_{2k}/M_{2k}^{(0)}$ decreases steadily towards $0$ as $k \to \infty$. For any fixed $k$, $M_{2k}/M_{2k}^{(0)} \to 1$ as $u \to \infty$.*

**Proof of Lemma 1 – first part.** For $k > 1$, $\frac{M_{2k}}{(2k)!} \cdot \frac{B_1^{k+1}}{(k+1)!} > \frac{M_{2(k+1)}}{(2(k+1))!} \cdot \frac{B_1^k}{k!}$, since each side is the sum of products of coefficients in (4.2) and (4.3) with the same numerators, but some of those on the right have larger denominators. Dividing by $\frac{B_1^{2k+1}}{k!(k+1)!}$ gives

$$\frac{M_{2k}}{M_{2k}^{(0)}} > \frac{M_{2(k+1)}}{M_{2(k+1)}^{(0)}} \tag{4.4}$$

Thus steady decrease of $M_{2k}/M_{2k}^{(0)}$ with $k$ from its values 1 at $k = 0$ and 1 is established. Decrease towards zero is confirmed in Section 5. Next, examination of some coefficients in (4.2) gives:

$$\frac{M_4}{M_4^{(0)}} = 1 - \frac{R_2}{2}; \quad \frac{M_6}{M_6^{(0)}} = 1 - \frac{3R_2}{2} + \frac{2R_3}{3}; \quad \frac{M_8}{M_8^{(0)}} = 1 - 3R_2 + \frac{3R_2^2}{4} + \frac{8R_3}{3} - \frac{11R_4}{8} \tag{4.5}$$

(see Leboeuf 2003 for further discussion in the case $q = 1$). Differences from 1 for larger $k$ are more complicated but finite expressions involving products of powers of the ratios $R_k$. Thus ratios indicate the difference of $P_u(v)$ from normal form for small values of $v$. Theorems 1 and 2 establish their behaviour, and confirm that as $u \to \infty$, all $R_{2k} \to 0$, and hence for any fixed $k > 1$ $M_{2k}/M_{2k}^{(0)} \to 1$, completing this part of the proof.

Indeed, (2.13) together with Theorem 2 indicates that $R_2$ is the fundamental determinant of difference of $P_u(v)$ from normal at small numbers of standard deviations, since expressions such as (4.5) are close to polynomials in $R_2$. This indicates that for an inter-residue race distribution determined by several characters, $P_u(v)$ will be closer to normal form than is a distribution determined by a single character – for example, that of a square/non square race – since $R_2$ is smaller.

The proof of Theorem 1, and other proofs later, uses the following result, clearly true by considering the meaning of the statements made:

**Lemma 2.** If $\xi(U) = O(\log U)$ then $\int_U^\infty \frac{\xi(w)}{w^{m+1}} dw = O\left(\frac{\log U}{mU^m}\right)$ (4.6)

**Proof of Theorem 1.** This draws on Lehmann 1966 Lemmas 1 and 2 and other work. The result

$$N_\chi(u) = \frac{u}{2\pi}\log\frac{q^*u}{2\pi} - \frac{u}{2\pi} + N_R(u) \quad = \frac{u(y-1)}{2\pi} + N_R(u) \quad (4.7)$$

where $N_R(u)$ is $O(\log u)$ and is a "saw tooth" function with upward discontinuities of 1 at $u = \gamma$ whose average is bounded (near 7/8 for $q = 1$: Ingham 1932 Theorem 25), has been proved for real characters, and also in terms of the average over the complex conjugate pairs which appear in (1.1). Using Stieltjes integrals and integrating by parts twice gives for $k \geq 1$:

$$b_k = \int_u^\infty \frac{N'_C(w)}{(1/4+w^2)^{2k}} dw = \frac{1}{2\pi}\int_u^\infty \frac{\log(q^*w/2\pi)}{w^{2k}} dw + \int_u^\infty \frac{N_R'(w)}{w^{2k}} dw + \varepsilon_k$$

$$= \frac{y + 1/(2k-1)}{2(2k-1)\pi u^{2k-1}} - \frac{N_R(u)}{u^{2k}} - \frac{2kN_I(u)}{u^{2k+1}} + 2k(2k+1)\int_u^\infty \frac{N_I(w)}{w^{2k+2}} dw + \varepsilon_k. \quad (4.8)$$

where $N_I(u) = \int_{u_0}^u N_R(w)dw$, $u_0$ being some value less than $u$ at which $N_I(u_0) = 0$. When $q = 1$, assuming the Riemann Hypothesis, Littlewood 1924 proved that $N_I(u) = O(\log u)$. Assuming this is generally true, the $N_I$ term in (4.8) is $O(ku^{-(2k+1)}\log u)$, as is the integral term (using (4.6)). The term $\varepsilon_k$ arises from the replacement of discrete values $(1/4 + \gamma^2)^{-2k}$ by continuous $w^{-2k}$. As the spacing of zeros is of order $2\pi/\log w$, the error introduced by this replacement is $O(4\pi k/(w^{2k+1}\log w))$, and $\varepsilon_k = O(u^{-2k}/\log u)$. Thus in (4.8) the $N_R$ term dominates other error terms. This demonstrates (2.11). (2.12) and (2.13) then follow algebraically.

The author has not seen results on the behaviour of $N_I(u)$ for a general $L$-function, but computations suggest that (2.11) applies generally. Indeed, the following argument suggesting higher accuracy in (2.11) when $q = 1$ is supported by statistical analysis for a variety of $L$-functions. As discussed by Odlyzko 1987, ultimately values of $N_R$ in the neighbourhood of $u = U$ are normally distributed, standard deviation $\sigma_E = \sqrt{\log\log U/2\pi^2}$. This quantity increases very slowly, being of order 0.3 for the largest values of $u$ considered in this paper. The observed root mean square (rms) deviation $\sigma_I$ of $N_I$ also appears to increase very slowly. Because $N_I(u)$ is oscillatory, whatever its ultimate behaviour the integral term in (4.8) should be no larger than the $N_I(U)$ term provided $2k \leq U$. So the proportionate rms errors in $b_k$ as given by (2.11) from $N_R$ and $N_I$ will be smaller than the bounds, being of order $\frac{2\pi(2k-1)\sigma_E}{uy}$ and $\frac{4\pi k(2k-1)\sigma_I}{u^2 y}$ respectively. Estimates of the rms error in ratios $r_k/r_2^{k-1}$ as given by (2.13) combine errors in the relevant $b_k$, noting that numerators and denominators move in parallel, to give $\frac{8\pi(k-1)(k-2)\sigma_E}{3(2k-1)uy^3}$ and $\frac{8\pi(k-1)(k-2)\sigma_I}{u^2 y}$ respectively. Since from (4.7) $uy \sim 2\pi N_\chi$, the proportionate rms error decreases faster than $1/N_\chi$, and increases with $k$ no faster than $k^2$. Note that the rms error in the estimate of $r_2$ by (2.12) is dominated by $\frac{4\pi\sigma_E}{uy}$ and decreases only as $1/N_\chi$.

To obtain (2.14), set $\gamma_0 <$ all $\gamma$, so that $N_\chi(\gamma_0) = 0$. Then using (4.7)

$$S_\chi(u) = 2 \int_{\gamma_0}^{u} \frac{N'_\chi(w)}{\sqrt{1/4 + w^2}} dw = \frac{\log u}{\pi} \left[ \frac{\log u}{2} + \log\left(\frac{q^*}{2\pi}\right) \right] - S_0 + S_1 + S_2 \quad \text{where}$$

$$S_0 = \frac{\log u_0}{\pi} \left[ \frac{\log u_0}{2} + \log\left(\frac{q^*}{2\pi}\right) \right]$$

$$S_1 = 2 \int_{\gamma_0}^{u} N'_\chi(w) \left( \frac{1}{\sqrt{1/4 + w^2}} - \frac{1}{w} \right) dw = 2N_\chi(u) \left( \frac{1}{\sqrt{1/4 + u^2}} - \frac{1}{u} \right)$$

$$-2 \int_{\gamma_0}^{\infty} N_\chi(w) \frac{d}{dw}\left( \frac{1}{\sqrt{1/4 + w^2}} - \frac{1}{w} \right) dw + 2 \int_{u}^{\infty} N_\chi(w) \frac{d}{dw}\left( \frac{1}{\sqrt{1/4 + w^2}} - \frac{1}{w} \right) dw$$

$N_\chi(u) = O(u \log u)$, $0 > \frac{1}{\sqrt{1/4+u^2}} - \frac{1}{u} > \frac{-1}{8} u^{-3}$ and $0 < \frac{d}{dw}\left[\frac{1}{\sqrt{1/4+w^2}} - \frac{1}{w}\right] < \frac{3}{8} w^{-4}$. Therefore, the above integrals exist, and the second integral is $O(u^{-2} \log u)$ using (4.6) Also

$$S_2 = 2 \int_{\gamma_0}^{u} \frac{N'_R(w)}{w} dw = 2 \left[ \frac{N_R(u)}{u} - \frac{N_R(\gamma_0)}{\gamma_0} + 2 \int_{\gamma_0}^{\infty} \frac{N_R(w)}{w^2} dw - 2 \int_{u}^{\infty} \frac{N_R(w)}{w^2} dw \right]$$

the variable terms of which are $O(u^{-1} \log u)$ and $O(u^{-2} \log u)$ respectively, again using (4.6). (2.14) follows, completing the proof. Computations give $\Delta_\chi \cong 0.50309$ for the prime count distribution, and -0.0836, -0.1224, and -0.2103 respectively for square/non-square races with $q = 4, 7$ and 13. The rms deviation of fluctuations about these levels appears to diminish more quickly than (2.14) suggests, as would follow from the above discussion of fluctuations in $N_R$.

**Proof of Theorem 2.** From (2.11), for larger $u$ the $b_k(\chi, u)$ for a particular distribution vary with $\chi$ only through differences in $y = \log\left(\frac{q^* u}{2\pi}\right)$ (and thus when $u$ is prime only through differences in the error terms in (2.11).) The results follow using (3.3).

## 5. BEHAVIOUR OF $P_u(v)$ FOR LARGE $v$

**Lemma 3: bounds on $I$ Bessel functions.** *For $w > 0$*

$$0 < (\log I_0(w))' < \min(w/2, 1), \qquad 0 < \log I_0(w) < \min(w^2/4, w) \qquad (5.1)$$

**Proof.** Define $j_l$ as the zeros of the Bessel function $J_0(x)$; $j_1 \cong 2.2048$, $j_2 \cong 5.5201$, etc.. Then

$$I_0(w) = \prod_{l=1}^{\infty}\left(1 + \frac{w^2}{j_l^2}\right); \qquad \log I_0(w) = \sum_{l=1}^{\infty} \log\left(1 + \frac{w^2}{j_l^2}\right) \qquad (5.2)$$

Therefore for $w > 0$, $(\log I_0(w))' = 2w \sum_{l=1}^{\infty} \frac{1}{j_l^2 + w^2}$. So $0 < (\log I_0(w))' < \frac{w}{2}$, since $\sum_{l=1}^{\infty} \frac{1}{j_l^2} = \frac{1}{4}$ by comparing the $w^2$ term in the product (5.2) and in the power series

$$I_0(w) = \sum_{k=0}^{\infty} \frac{w^{2k}}{2^{2k}(k!)^2} \qquad (5.3)$$

$(\log I_0(w))' = I_1(w)/I_0(w) < 1$ since $I_0(w) - I_1(w) = \frac{1}{\pi}\int_0^{\pi}(1 - \cos\theta)\exp(w\cos\theta)\,d\theta > 0$.
The second result follows by integration. From (2.8) follows immediately that $L_u(s) < \frac{1}{2}s^2\sigma_u^2$.

**Proof of Theorem 3.** Following Montgomery 1979, for <u>any</u> positive $v$ and $s$, with $D_0(z)$ defined by (3.1)

$$E(v) = \lim_{Z \to \infty} \frac{1}{Z}\int_0^Z dz: D_0(z) > v \qquad \leq e^{-sv} \lim_{Z \to \infty} \frac{1}{Z}\int_0^Z \exp(D_0(z))\,dz: D_0(z) > v$$

$$< e^{-sv} \lim_{Z \to \infty} \frac{1}{Z}\int_0^Z \exp(D_0(z))\,dz \qquad = e^{-sv}\bar{P}_0(s)$$

so that $\log E(v) < L_0(s) - sv$. Also for <u>any</u> $u > 0$ and not equal to any $\gamma$

$$L_0(s) = \sum_\chi \sum_{\gamma < u} \log I_0\left(\frac{2\alpha_\chi s}{\sqrt{1/4 + \gamma^2}}\right) + L_u(s) < sS(u) + \frac{1}{2}s^2\sigma_u^2 \qquad (5.4)$$

so that

$$\log E(v) < \frac{1}{2}s^2\sigma_u^2 - s(v - S(u)) = \frac{\sigma_u^2}{2}\left(s - \frac{v - S(u)}{\sigma_u^2}\right)^2 - \frac{(v - S(u))^2}{2\sigma_u^2}$$

(2.15) follows, choosing $s = \frac{v - S(u)}{\sigma_u^2}$, permissible if $v \geq S(u)$. Montgomery takes $v = 2S(u)$, so that $\log E(2S(u)) < -\frac{(S(u))^2}{2\sigma_u^2}$. Taking $v = 3S(u)$ gives a generally better bound:

$$\log E(3S(u)) < -\frac{2(S(u))^2}{\sigma_u^2} . \qquad (5.5)$$

Whatever is chosen, there is a series of values of $v$, from (2.14) increasing with $u$ like $(\log u)^2$,

for which using (2.11) and (2.14) $\log E(v)$ behaves like $-u(\log u)^3$. This decrease is faster than would be given by any normal distribution. However, it is often convenient to set $u = 0$, giving

$$\log E(v) < -\frac{v^2}{2\sigma_0^2} \tag{5.6}$$

as quoted by Montgomery and Odzylko 1988 (Equation (5)).

**Proof of Lemma 1: second part.** From (4.4), $M_{2k}/M_{2k}^{(0)} \to \varepsilon$ from above as $k \to \infty$, where $0 \leq \varepsilon < 1$. Therefore $\bar{P}_u(s) > \varepsilon \sum_{k=0}^{\infty} \frac{M_{2k}^{(0)} s^{2k}}{(2k)!} = \varepsilon \exp(\sigma_u^2 s^2)$. If $\varepsilon > 0$, $\frac{L_u(s)}{s^2} > \frac{\log \varepsilon}{s^2} + \sigma_u^2 > \sigma_u^2/2$ for large enough $s$. However, similarly to (5.4), for any $u' > u$, $L_u(s) = \sum_\chi \sum_{u < \gamma < u'} \log I_0\left(\frac{2\alpha_\chi s}{\sqrt{1/4+\gamma^2}}\right) + L_{u'}(s) < s(S(u') - S(u)) + \frac{1}{2}s^2 \sigma_{u'}^2$. From (2.11) and (2.14), choosing $u'$ proportional to $s$, $L_u(s) = O(s(\log s)^2)$. Therefore, $L_u(s)/s^2 \to 0$ as $s \to \infty$. So $\varepsilon = 0$ and $M_{2k}/M_{2k}^{(0)} \to 0$ as $k \to \infty$.

## 6. PRECISE CALCULATIONS OF $P_0(v)$ FOR SMALL $v$

**Expansion of $\log J_0(w)$.** Using some results given in Fiorilli and Martin 2012 Lemma 2.8, since $J_0(w) = 1 - \sum_{n=1}^{\infty} \frac{(-1)^{n-1} w^{2n}}{2^{2n}(n!)^2}$, is absolutely convergent, $\log(J_0(w))$ may be expanded as a power series in $w$ for the range $w < j_1$ and being an even function, zero at $w = 0$, is of the form

$$\log J_0(z) = -\sum_{k=1}^{\infty} \frac{c_k w^{2k}}{2^k} . \tag{6.1}$$

Some algebra gives the values of $c_k$ indicated in the statement of Theorem 4. (In terms of Fiorilli and Martin's notation, $c_k = -2^k \lambda_{2k}$.) Since $J_0(w) = \prod_{l=1}^{\infty} \left(1 - \frac{w^2}{j_l^2}\right)$, taking logarithms and expanding,

$$c_k = \frac{2^k}{k} \sum_{l=1}^{\infty} j_l^{-2k} > 0. \tag{6.2}$$

Examination of tables of $j_l$, and for larger $l$ the asymptotic formula (9.5.12) in Abramowitz and Stegun 1964, shows that the difference $j_{l+1} - j_l$ increases steadily towards $\pi$ and exceeds 3.13 for $l \geq 2$. Hence for $l > 2$, $j_l > j_2 + 3.13(l-2) > j_2(1 + 0.567(l-2))$. So:

$$c_k - \frac{1}{kj_1^{2k}} = \frac{1}{kj_2^{2k}}\left[1 + \sum_{l>2}\left(\frac{j_2}{j_l}\right)^{2k}\right] < \frac{1}{kj_2^{2k}}\left[1 + \sum_{l>2}\left(\frac{1}{1 + 0.567(l-2)}\right)^{2k}\right]$$

$$< \frac{1}{kj_2^{2k}}\left[1 + \int_0^{\infty}\left(\frac{1}{1 + 0.567x}\right)^{2k} dx\right] = \frac{1}{kj_2^{2k}}\left[1 + \frac{1}{0.567(2k-1)}\right] < \frac{1.16}{kj_2^{2k}}$$

for $k > 5$, which gives the good approximation (2.17), noting that $j_1^2/j_2^2 < 0.19$. The ratio $c_k/c_{k+1}$ tends from above to $j_1^2/2 \cong 3.134$.

**Proof of Theorem 4.** When (6.1) is valid for all terms in the exponent:

$$\hat{P}_u(\omega) = \exp\left(\sum_\chi \sum_{\gamma>u} \log J_0\left(\frac{2\alpha_\chi \omega}{\sqrt{1/4+\gamma^2}}\right)\right) = \exp\left(-\sum_\chi \sum_{\gamma>u} \sum_{k=1}^\infty c_k \left(\frac{2\alpha_\chi^2 \omega^2}{1/4+\gamma^2}\right)^k\right)$$

$$= \exp\left(-\sum_{k=1}^\infty 2^k c_k B_k \omega^{2k}\right) = \exp\left(-\sum_{k=1}^\infty c_k R_k \tau^{2k}\right) \quad (6.3)$$

where $\tau = \sqrt{2B_1}\omega$, and as (2.10) implies, $R_1 = 1$. ((6.3) is equivalent to Proposition (2.12) of Fiorilli and Martin 2012.)

In the single series case ($R_k = r_k$), $\sum_{k=1}^\infty c_k R_k \tau^{2k}$ converges for $|\tau| < T$ because $c_{k+1}R_{k+1}/c_k R_k \to 1/T^2$. In the multiple series case, the approximate zone of convergence is established as (2.18) by Theorem 2. (6.3) shows how as $u \to \infty$, $\hat{G}_u(\omega)$ as a function of $\tau$ tends (though not uniformly) to the Fourier transform of a normal distribution. For $|\tau| < T$, the proportionate error in (2.16) is given by assuming $c_{k+1}R_{k+1}/c_k R_k = 1/T^2$ for $k \geq K$.

**Application to Rubinstein and Sarnak's method.** Accordingly

$$\hat{P}_0(\omega) \cong \prod_\chi \prod_{0<\gamma<u} J_0\left(\frac{2\alpha_\chi \omega}{\sqrt{1/4+\gamma^2}}\right) \exp\left(-\sum_{k=1}^K 2^k c_k B_k \omega^{2k}\right) \quad (6.4)$$

with an error that can be estimated from (2.16). The Poisson summation formula gives, for any positive real $\Delta\omega$:

$$\sum_{l=-\infty}^\infty P_0\left(\frac{2\pi l}{\Delta\omega}+v\right) = \frac{\Delta\omega}{2\pi}\sum_{m=-\infty}^\infty \hat{P}_0(m\Delta\omega)e^{imv\Delta\omega} = \frac{\Delta\omega}{2\pi}\sum_{m=-\infty}^\infty \hat{P}_0(m\Delta\omega)\cos(mv\Delta\omega) \quad (6.5)$$

since $\hat{P}_0(\omega)$ is symmetric. Integrating over $0$ to $v$, and using symmetry of $P_0(v)$

$$E(v) - \epsilon(v) = \frac{1}{2} - \frac{v\Delta\omega}{2\pi} - \frac{1}{\pi}\sum_{m=1}^\infty \hat{P}_0(m\Delta\omega)\frac{\sin(mv\Delta\omega)}{m} \quad (6.6)$$

where
$$\epsilon(v) = \sum_{l=1}^\infty \left[E\left(\frac{2\pi l}{\Delta\omega}-v\right) - E\left(\frac{2\pi l}{\Delta\omega}+v\right)\right] = \sum_{l=1}^\infty \int_{\frac{2\pi l}{\Delta\omega}-v}^{\frac{2\pi l}{\Delta\omega}+v} P_0(v')\,dv' \quad (6.7)$$

Provided $2\pi/\Delta\omega > v$, $0 < \epsilon(v) < E\left(\frac{2\pi}{\Delta\omega}-v\right)$ since the ranges of integration in (6.7) do not overlap. Therefore, selecting $\Delta\omega$ so that $\exp\left(\frac{-1}{2\sigma^2}\left(\frac{2\pi}{\Delta\omega}-v\right)^2\right)$ is negligible, and using (5.6), $\epsilon(v)$ can be left out of (6.6). The sum on the right hand side of (6.6) can be limited to $|m\Delta\omega| < C$, where $C$ is such that $\hat{P}_0(\omega)$ is negligible for $\omega > C$. Then:

$$E(v) \cong \frac{1}{2} - \frac{v\Delta\omega}{2\pi} - \frac{1}{\pi}\left(\sum_{0<m\Delta\omega<C} \hat{P}_0(m\Delta\omega)\frac{\sin(mv\Delta\omega)}{m}\right) \quad (6.8)$$

Similarly, from (6.5) provided $P_0(2\pi/\Delta\omega - v)$ is negligible

$$P_0(v) \cong \frac{\Delta\omega}{\pi}\left(\frac{1}{2} + \sum_{0<m\Delta\omega<C} \hat{P}_0(m\Delta\omega)\,\cos(mv\Delta\omega)\right) \qquad (6.9)$$

**Economy in computation.** Three factors allow calculations using (6.8) and (6.9) to be done much more economically than has been the case in previous applications of Rubinstein and Sarnak's method:

a) For small values of $K$, perhaps 5 or 7, (2.16) is already <u>far</u> more accurate than the forms of remainder previously employed, which in present notation amount to $\hat{P}_u(\omega) \cong 1 - B_1\omega^2$, or sometimes $\hat{P}_u(\omega) \cong 1 - B_1\omega^2 + (B_1^2/2 - B_2/4)\omega^4$ (See (4.15) of Rubinstein and Sarnak 1994). These are the initial terms in the power series expansion of (2.16), and take no account of the approximately normal nature of $P_u(v)$. (2.16) allows high accuracy with $u$ such that $N(u) \leq 10$, rather than many thousands;

b) (5.6) gives assurance that $\epsilon(v) < E(2\pi/\Delta\omega - v)$ can be neglected with much larger values of $\Delta\omega$ than have commonly been employed, thus much reducing the number of terms in (6.8), and similarly for (6.9)

c) As $\tau$ approaches T and the series for $\log \hat{P}_u(\omega)$ approaches divergence and (2.16) less accurate, its sum becomes large and negative, so $\hat{P}_u(\omega)$ become negligibly small. Therefore, a few terms of the series remain sufficient to give accurate results in (6.8) and (6.9). There could be contributions from values $|\tau| > T$ for which the representation (2.16) fails. However, for such values, let $u' > u$ be such that $\sum_{k=1}^{\infty} 2^k c_k B_k(u')\omega^{2k}$ just converges (which value exists since $R_2(u) \to 0$ as $u \to \infty$), and $N'$ be the number of terms with $u < \gamma < u'$. For each of these terms, $\frac{2\alpha_\chi\omega}{\sqrt{1/4+\gamma^2}} > j_1$ and therefore $\left|J_0\left(\frac{2\alpha_\chi\omega}{\sqrt{1/4+\gamma^2}}\right)\right| < 0.41$. Therefore $|\hat{P}_u(\omega)| < (0.41)^{N'}\hat{P}_{u'}(\omega)$, and $\hat{P}_{u'}(\omega)$ will itself be very small, so a contribution can be neglected.

All calculations described here have been repeated with different values of $N(u)$, $\Delta\omega$ and $C$, with results unchanged to the accuracy given.

**Computations of prime count distribution.** Table 3 sets out a *complete* calculation of $E(1)$, the limiting logarithmic density of $\pi(x) > \text{Li}(x) + \text{Li}(\sqrt{x})$, to 5 significant figures ($10^{-11}$), using $\Delta\omega = \pi/2$, so that half the terms in (6.8) are zero. Quantities are displayed to fewer digits than calculated. The maximum error estimated by (2.16) in individual terms of (6.8) is less than $10^{-12}$. The calculation depends on $E(3)$ being negligible; (5.6) estimates it as less than $4.9 \times 10^{-43}$, and (5.5) as less than $1.12 \times 10^{-62}$ (its actual value is near $1.603 \times 10^{-95}$ as given in Table 1). Just 60 calculations of Bessel functions are needed, whilst Rubinstein and Sarnak used about 120 million to obtain 2 significant figures ($u = 88{,}190$, $N(u) = 120{,}000$, $C = 50$, $\Delta\omega = 0.05$).

As $v$ increases, more terms are required to deliver a certain proportional accuracy. Figure 2 shows this in regard to sums (6.9) for $P_0(v)$ up to the value of $m$ indicated (in this case including both even and odd $m$).

TABLE 3. CALCULATION OF $E(1)$ WITH $u = 35$ $(N(u) = 5), K = 7, \Delta\omega = \pi/2$

| $m$ | $\omega$ | $\tau/T$ | $\hat{f}_{\gamma 1}$ | $\hat{f}_{\gamma 2}$ | $\hat{f}_{\gamma 3}$ | $\hat{f}_{\gamma 4}$ | $\hat{f}_{\gamma 5}$ | $\hat{F}_u$ | $\hat{P}_u$ | Sum for $E(1)$ |
|---|---|---|---|---|---|---|---|---|---|---|
|  |  |  |  |  |  |  |  |  |  | 0.25000000000 |
| 1 | 1.571 | 0.037 | 0.9877 | 0.9944 | 0.9961 | 0.9973 | 0.9977 | 0.9735 | 0.9703 | -0.05066067594 |
| 3 | 4.712 | 0.110 | 0.8920 | 0.9504 | 0.9648 | 0.9762 | 0.9796 | 0.7822 | 0.7618 | 0.01256921934 |
| 5 | 7.854 | 0.183 | 0.7146 | 0.8653 | 0.9038 | 0.9345 | 0.9439 | 0.4930 | 0.4690 | -0.00215114117 |
| 7 | 10.996 | 0.257 | 0.4810 | 0.7447 | 0.8159 | 0.8736 | 0.8916 | 0.2277 | 0.2257 | 0.00018592266 |
| 9 | 14.137 | 0.330 | 0.2244 | 0.5966 | 0.7052 | 0.7955 | 0.8241 | 0.0619 | 0.0845 | 0.00000085594 |
| 11 | 17.279 | 0.403 | -0.0198 | 0.4306 | 0.5769 | 0.7026 | 0.7432 | -0.0026 | 0.0245 | -0.00000096511 |
| 13 | 20.420 | 0.477 | -0.2196 | 0.2573 | 0.4368 | 0.5979 | 0.6511 | -0.0096 | 0.0054 | 0.00000031044 |
| 15 | 23.561 | 0.550 | -0.3511 | 0.0877 | 0.2913 | 0.4845 | 0.5502 | -0.0024 | 0.0009 | 0.00000026436 |
| 17 | 26.704 | 0.623 | -0.4021 | -0.0680 | 0.1468 | 0.3661 | 0.4432 | 0.0007 | 0.0001 | 0.00000026298 |
| 19 | 29.845 | 0.697 | -0.3737 | -0.2006 | 0.0098 | 0.2461 | 0.3329 | 0.0001 | 0.0000 | 0.00000026299 |
| 21 | 32.987 | 0.770 | -0.2791 | -0.3030 | -0.1142 | 0.1282 | 0.2222 | -0.0003 | 0.0000 | 0.00000026300 |
| 23 | 36.128 | 0.843 | -0.1414 | -0.3705 | -0.2199 | 0.0158 | 0.1139 | 0.0000 | 0.0000 | 0.00000026300 |

$\omega = m\pi/2$; $\hat{f}_\gamma = J_0\left(2\omega/\sqrt{1/4 + \gamma^2}\right)$; $E(1) = 1/2 - 1/4 - \sum_{\text{odd } m > 0}(-1)^{(m-1)/2}\hat{F}_u\hat{P}_u/m\pi$

Including 3 more terms in the sum over $m$ than Table 1 shows takes contributions to (6.8) below $10^{-20}$. Values of $E(1)$ within $10^{-16}$, the approximate accuracy of computation, of $2.629967324 \times 10^{-7}$ were obtained with various values of $\Delta\omega$ and with $N(u) = 10, K = 10$; $N(u) = 25, K = 7$; $N(u) = 50, K = 5$; and $N(u) = 100, K = 4$. (For $k > 5$, $c_k$ were calculated with a few terms of $\frac{1}{k}\sum_{l=1}^{\infty} 2^k j_l^{-2k}$.) This "trade off" in accuracy between $N(u)$ and $K$ reflects the fact that with more initial terms, lower values of $\tau/T$ are used, and hence (2.16) converges more quickly. Figure 3 shows this in more detail.

FIGURE 2. % DIFFERENCES FROM $P_0(v)$ OF SUM TO $m$ TERMS, $K = 7, \Delta\omega = \pi/2$

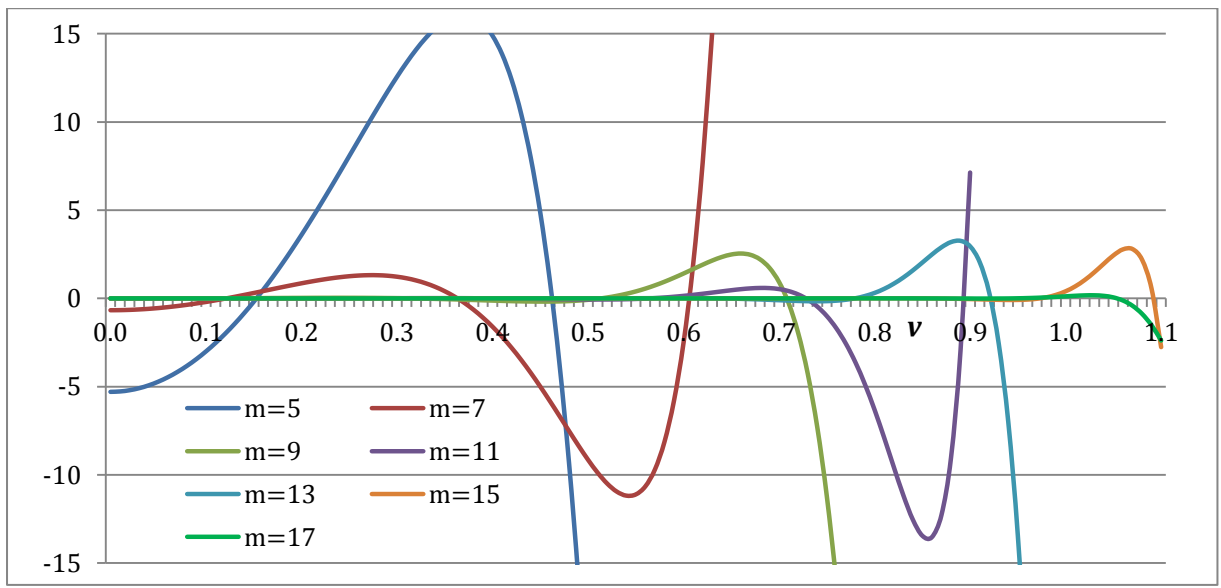

Data used are discussed below: changes within its known accuracy have effects no larger than of order $10^{-16}$. It can be stated with some confidence that for the prime count distribution:

$$2.62996732 \times 10^{-7} < E(1) < 2.62996733 \times 10^{-7} \tag{6.10}$$

FIGURE 3. VARIATION OF ERROR IN $E(1)$ WITH $N(u)$ AND $K$

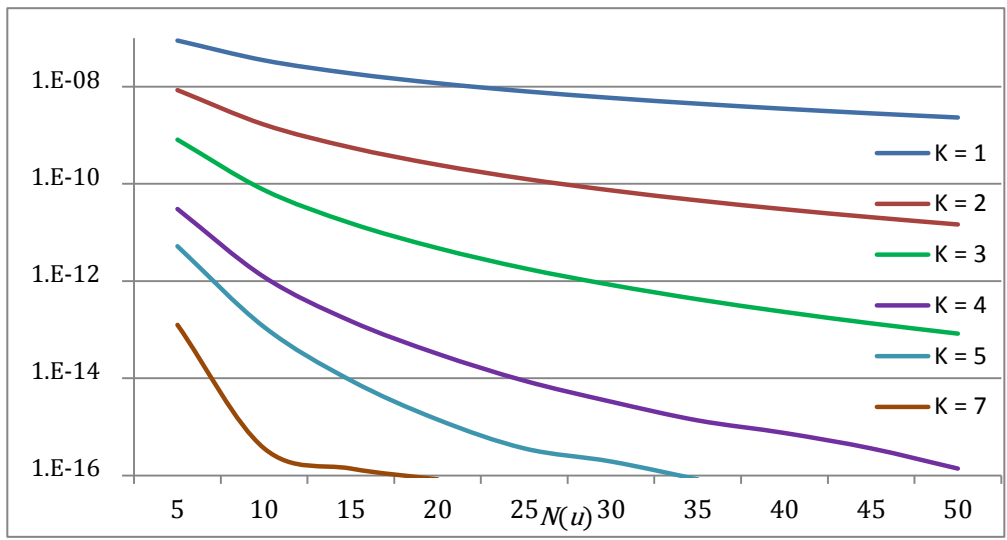

**Prime number races.** The results in Table 2 for races with $q$ prime or equal to 4 were calculated with $\Delta\omega = 1/2\sigma_0$. This requires $E(4\pi\sigma_0 - 1)$ to be negligible; from (5.6) its magnitude is less than $\exp\left(-\frac{1}{2}\left(4\pi - \frac{1}{\sigma_0}\right)^2\right) < 10^{-31}$ for all $\sigma_0$ quoted. Results to 7 decimal places can be obtained by summing about 12 terms in (6.8) with the values of $N_\chi(u)$ and $K$ given.

The ability to obtain accurate inter-residue race results with $N_\chi(u) = 1$ reflects the fact that the combination of series makes these distributions closer to normal, as discussed above and indicated by the smaller values of $\beta$ for these races by comparison with square/non-square races. Although first zeros of each $L$-function have high variability, even the second zeros show much less variability; for $u$ exceeding all first zeros, $P_u$ is close to normal and varies little amongst races. For example, with $q = 13$ a first zero is near 0.88, giving $\beta \cong 0.80$ for this component series. $\alpha_\chi$ for this zero is near 0.07, 0.5 and 0.93 in the 3 races and this variation accounts for almost all the differences in $\sigma_0$ and $\beta$ shown. However, if $u$ just exceeds all first zeros, for the 3 races $\sigma_u$ remains between 0.8 and 0.81 and $R_2(u)$ between 0.0118 and 0.0131.

Two cases for $q = 8$ are also given: the race between 1 and 3 requires calculation of $E(2)$, whilst the method for the three way race follows that set out in Feuerverger and Martin 2000. Good results are obtained with small $N$ and $K$, and $\Delta\omega = 0.8$.

Results can be obtained to within the accuracy of computation or data, for example 0.004072076720775 for the modulo 4 race with $N = 50, K = 5$.

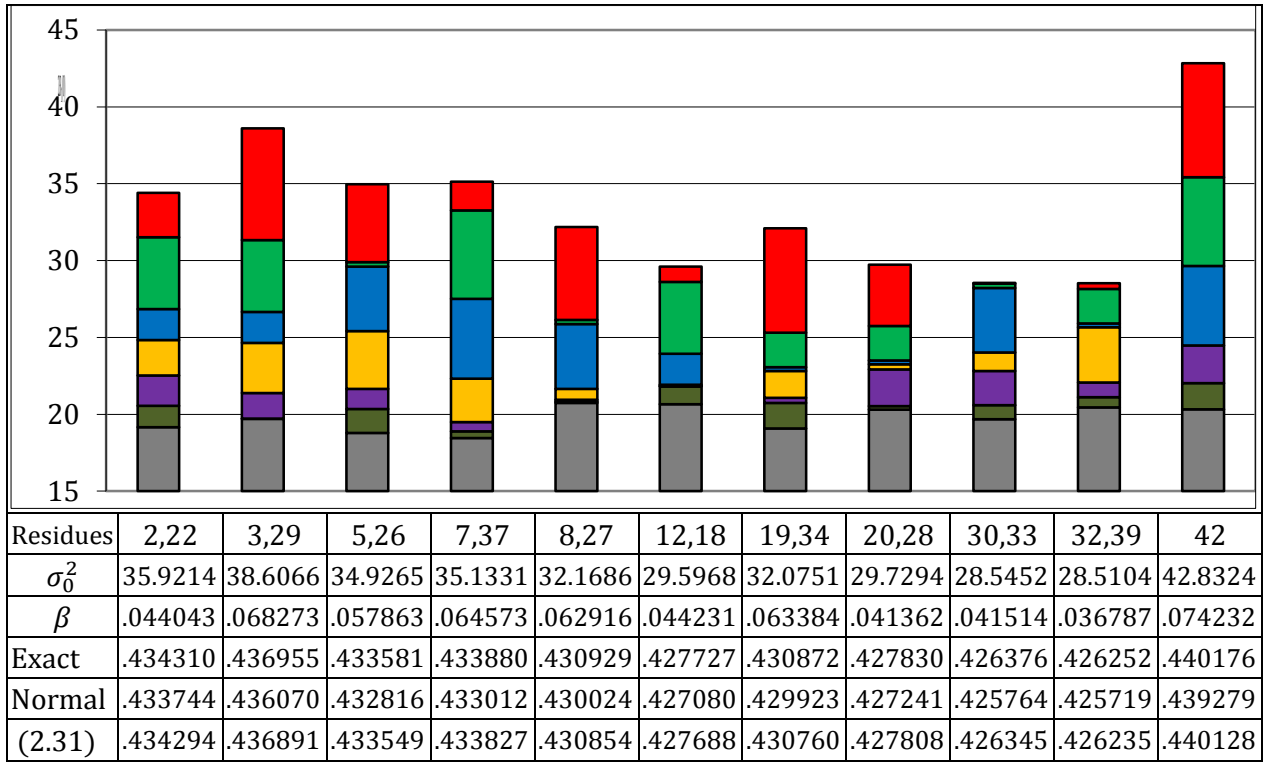

FIGURE 4. MODULO 43 – RACES 1 VS. OTHER RESIDUES

| Residues | 2,22 | 3,29 | 5,26 | 7,37 | 8,27 | 12,18 | 19,34 | 20,28 | 30,33 | 32,39 | 42 |
|---|---|---|---|---|---|---|---|---|---|---|---|
| $\sigma_0^2$ | 35.9214 | 38.6066 | 34.9265 | 35.1331 | 32.1686 | 29.5968 | 32.0751 | 29.7294 | 28.5452 | 28.5104 | 42.8324 |
| $\beta$ | .044043 | .068273 | .057863 | .064573 | .062916 | .044231 | .063384 | .041362 | .041514 | .036787 | .074232 |
| Exact | .434310 | .436955 | .433581 | .433880 | .430929 | .427727 | .430872 | .427830 | .426376 | .426252 | .440176 |
| Normal | .433744 | .436070 | .432816 | .433012 | .430024 | .427080 | .429923 | .427241 | .425764 | .425719 | .439279 |
| (2.31) | .434294 | .436891 | .433549 | .433827 | .430854 | .427688 | .430760 | .427808 | .426345 | .426235 | .440128 |

**Modulo 43 races.** As $q$ increases further, it is not even necessary to take explicit account of all initial zeros, only of the smallest. With $q = 43$, race calculations can involve up to 41 Dirichlet characters, but results to 6 decimal places, as set out below Figure 4 (row "Exact"), can be obtained with $K = 3$ by explicitly representing only the 7 first zeros with imaginary parts of magnitude less than 1 (the smallest being about 0.138). For these non-extreme races, almost all of the variation in result with residue class arises from changes in standard deviation; the correlation between the exact result and its value for a normal distribution with the same variance $\sigma_0^2$ (row "Normal") is 0.9996. Figure 4 shows how almost all variation amongst races in $\sigma_0^2$ (and in $\beta$ which has correlation 0.79 with $\sigma_0^2$) is driven by the varying $\alpha_\chi$ for 6 first zeros, represented by red, green, blue, orange, violet and olive. Such dependence has been commented on by Bays *et al*. 2001. (The other first zero explicitly represented is of the *L*-function determining the square/non-square distribution. This distribution has $\beta \cong 0.69$, and its density is so far from normal as to be "double humped"). Results from (2.31) are also given, and are discussed in Section 10.

**Calculation of $B_k$ and $R_k$.** The most accurate values of zeros $\gamma$ available were taken from tables on the internet (Odlyzko, de Silva). For $k > 1$, $b_k$ was evaluated by summation of powers of a convenient number of $\gamma$ (1000000 for zeta zeros, 10000 for others) and remaining contributions estimated using (2.11). Where an accurate value of $B_1(0)$ was available, (eg $1 + \gamma/2 - \log(4\pi)/2$ for the zeta zeros, where $\gamma$ here denotes Euler's constant, and in other cases as given in Rubinstein and Sarnak 1994), $b_1(u)$ was calculated by subtraction from it. $b_1$ was also estimated by summation as for $k > 1$, with a difference within the estimated error in (2.11) when both results were available. The key driver of data error in calculations is error in $b_1$. It was checked that adjusting data within its error range had no effect on results within quoted accuracy. The calculations for $q = 43$ used the first 250 zeros of each $L$-function, calculated using Rubinstein's zero finder. Similar checks on accuracy were made. Bessel functions needed to be calculated for arguments up to about 8; for these, the power series was used and checked against tables as giving 15 decimal places.

These calculations reflect availability of sufficient values of $\gamma$ to obtain high accuracy in $b_k$. However, all $b_k(0)$ can in principle be calculated without reference to these values (see Fiorilli and Martin 2012 Lemma 3.15) and other $b_k(u)$ obtained by subtraction. Given the rapid decrease of $(1/4 + \gamma^2)^{-k}$ for larger $k$, this method would be prone to inaccuracy; but its existence demonstrates that the calculations here do not depend explicitly on individual $\gamma > u$.

## 7. METHOD OF STEEPEST DESCENT: BEHAVIOUR OF LAPLACE TRANSFORM

**Basis of the method.** The Rubinstein-Sarnak method in the form above can give accurate results for highly biased races (see Fiorilli 2012). However, because it obtains the difference of quantities of order 1, its absolute accuracy is limited by that of computation, here about $10^{-16}$. For more extreme results, a different method is needed.

(5.5) suggests that $\log E(v)$ (and hence $\log P_0(v)$) eventually decrease at an accelerating rate. These are the conditions suited to using the method of steepest descent in the form originally introduced by Laplace. For any $s$, the product $P_0(v)e^{sv}$ will rise to a maximum, then decline; the value of $\bar{P}_0(s)$ will be largely determined by the values of $P_0(v)$ at and near this maximum. So it should be possible to estimate $P_0(v)$ from $\bar{P}_0(s)$. The principle of the method is set out below, following which the validity of approximations made is considered. Define

$$Q(v) = \log P_0(v), \quad V(s) = L'_0(s) = \frac{\bar{P}_0'(s)}{\bar{P}_0(s)}, \quad \eta(s) = 1/\sqrt{-Q''(V)}, \quad \Delta s = s + Q'(V) \qquad (7.1)$$

so that accelerating decline means $Q''(v) < 0$. Writing $v' = v - V$,

$$\int_{-\infty}^{\infty} v' P_0(v) e^{sv} dv = \int_{-\infty}^{\infty} v P_0(v) e^{sv} dv - V \bar{P}_0(s) = 0$$

so that $V(s)$ is the mean of the distribution $P_0(v)e^{sv}$. Also

$$P_0(v)e^{sv'} = \exp\left[Q(V) + v'\Delta s + \sum_{n=2}^{\infty} \frac{v'^n}{n!} Q^{(n)}(V)\right] \cong \left[1 + v'\Delta s + \frac{v'^3}{6} Q'''(V)\right] \exp\left[Q(V) - \frac{v'^2}{2\eta^2}\right]$$

where terms above the third in the Taylor expansion have been neglected and the exponentials of $v'$ and $v'^3$ terms approximated. For this to be valid, $\Delta s$ (which is related to the difference between the mean and the mode of the of the distribution of $P_0(v)e^{sv}$) must be small compared to the standard deviation $\eta$, and $\kappa = \eta^3 Q'''(V)$, (which is the skewness of this distribution) must also be small. Then

$$\bar{P}_0(s) \cong \frac{e^{sV}}{\pi} \int_{-\infty}^{\infty} \left[1 + v'\Delta s + \frac{v'^3}{6} Q'''(V)\right] \exp\left[Q(V) - \frac{v'^2}{2\eta^2}\right] dv' = \frac{P_0(V)\eta e^{sV}}{\sqrt{2\pi}}$$

$$\bar{P}_0'(s) \cong \frac{e^{sV}}{\pi} \int_{-\infty}^{\infty} [V + v'] \left[1 + v'\Delta s + \frac{v'^3}{6} Q'''(V)\right] \exp\left[Q(V) - \frac{v'^2}{2\eta^2}\right] dv'$$

$$= \bar{P}_0(s)\left[V + \frac{\eta}{\sqrt{2\pi}}\left(\eta\Delta s + \frac{\kappa}{2}\right)\right]$$

So $\Delta s \cong -\kappa/2\eta$, confirming consistency of the two assumptions of smallness. Also

$$Q''(V) = \frac{d}{dV}(\Delta s - s) = \frac{d}{dV}\left(\frac{-Q'''(V)}{2(Q''(V))^2}\right) - \frac{1}{L_0''(s)} = \kappa^2 - \frac{Q''''(V)}{2(Q''(V))^2} - \frac{1}{L_0''(s)}$$

The second term can be omitted under the approximations already used. So to order $\kappa$, $\eta = \sqrt{L_0''(s)}$ and

$$Q(V) \cong L_0(s) - sV - \frac{1}{2}\log(2\pi L_0''(s)); \quad P_0(V) = \bar{P}_0(s)e^{-sV}/\sqrt{2\pi L_0''(s)} \qquad (7.2)$$

If therefore accurate estimates of $L_0(s)$ and its derivatives can be made, $V(s)$ and the corresponding $P_0(V)$ may be calculated for any $s$.

To validate the approximations made in deriving (7.2) requires a fuller calculation to higher order, effectively estimating of the magnitude of further terms in the asymptotic expansion which gives solutions in the method of steepest descent. Theorem 5 and the estimates (2.19) are used to support this calculation. Also needed is confirmation for any $V \geq 0$, there exists a unique $s$ such that $L_0'(s) = V$, and that $L_0''(s)$, which estimates the variance of the near normal distribution of $P(v)e^{sv}$, is positive. This confirmation follows from (2.8) and

**Lemma 4.** $(\log(I_0(w)))'$ increases steadily from $0$ towards $1$ in $0 \leq w < \infty$, and
$$0 < (\log(I_0(w)))'' < 1/2 \qquad (7.3)$$

**Proof.** If $a > b \geq 0$,

$$I_1(a)I_0(b) - I_0(a)I_1(b) = \frac{1}{\pi^2} \int_{-1}^{1}\int_{-1}^{1} \frac{(w_1 - w_2)\exp(aw_1 + bw_2)}{\sqrt{(1-w_1^2)(1-w_2^2)}} dw_1 dw_2$$

$$= \frac{1}{\pi^2} \int_{-1}^{1}\int_{-1}^{1} \frac{(w_1 - w_2)\exp\frac{1}{2}[(a-b)(w_1 - w_2) + (a+b)(w_1 + w_2)]}{\sqrt{(1-w_1^2)(1-w_2^2)}} dw_1 dw_2$$

Averaging with the result of interchanging $w_1$ and $w_2$ gives

$$\frac{1}{\pi^2} \int_{-1}^{1}\int_{-1}^{1} \frac{(w_1 - w_2)\exp[\frac{1}{2}(a+b)(w_1 + w_2)]\sinh[\frac{1}{2}(a-b)(w_1 - w_2)]}{\sqrt{(1-w_1^2)(1-w_2^2)}} dw_1 dw_2 > 0$$

since the integrand is always positive. So $\frac{I_1(a)}{I_0(a)} > \frac{I_1(b)}{I_0(b)}$, $(\log(I_0(w)))'$ is an increasing function of $w$ in $w \geq 0$, and $(\log(I_0(w)))'' > 0$. From (5.1) $(\log(I_0(w)))' \to 1 - \varepsilon$ as $w \to \infty$ where $0 \leq \varepsilon < 1$. Since $I_1(0) = 0$, $(\log(I_0(w)))' = 0$ at $w = 0$. So, integrating, $\log(I_0(w)) - w < -\varepsilon w$. But

$$e^{-w}I_0(w) = \frac{e^{-w}}{\pi}\int_0^\pi \exp(w\cos\theta)\,d\theta = \frac{1}{\pi}\int_0^\pi \exp\left(-2w\sin^2\left(\frac{\theta}{2}\right)\right)d\theta > \frac{1}{\pi}\int_0^\pi e^{-w\theta^2/2}\,d\theta$$

$$> \frac{1}{\pi}\int_0^\infty e^{-w\theta^2/2}\,d\theta - \frac{1}{\pi}\int_\pi^\infty e^{-\pi w\theta/2}\,d\theta = \frac{1}{\sqrt{2\pi w}} - \frac{2}{\pi^2 w} > \frac{1}{2\sqrt{2\pi w}}$$

for $w > 32/\pi^3 \cong 1.032$. So then $\log(I_0(w)) - w > -\frac{1}{2}\log w - \log(8\pi)$, possible only if $\varepsilon = 0$. Also from (5.2)

$$(\log I_0(w))'' = 2\sum_{l=1}^\infty \frac{1}{j_l^2 + w^2}\left(1 - \frac{2w^2}{j_l^2 + w^2}\right) \leq 2\sum_{l=1}^\infty \frac{1}{j_l^2} = \frac{1}{2}$$

concluding the proof.

From Lemma 4 follows results analogous to (5.4); for any $u \geq 0$

$$0 < L_0'(s) \leq s\sigma_u^2 + S(u);\quad 0 < L_u''(s) \leq \sigma_u^2 \tag{7.4}$$

Theorem 5 also requires

**Lemma 5.** *For $|w| \leq 1$*

$$\log I_0(w) = O(w^2);\quad \frac{d^{2M}}{dw^{2M}}\log I_0(w) = O(1);\quad \frac{d^{2M-1}}{dw^{2M-1}}\log I_0(w) = O(w)\ (M \geq 1) \tag{7.5}$$

*Defining $I(w) = \log I_0(w) - w$, for large $w$*

$$I(w) = O(\log w);\quad \frac{d^M I(w)}{dw^M} = \frac{(-1)^M(M-1)!}{2w^M} + O(w^{-(M+1)}) \tag{7.6}$$

**Proof.** (7.5) follows from absolute convergence of $\log I_0(w) = \sum_{k=1}^\infty (-1)^{k-1}c_k w^{2k}/2^k$ when $|w| \leq 1$. That $I(w) = O(\log w)$ has been demonstrated in the proof of Lemma 4. For $M = 1$ (7.6) follows from permissible division of the asymptotic expansions

$$I_\nu(w) = \frac{e^w}{\sqrt{2\pi w}}\left(1 + \cdots + \frac{(1^2 - 4\nu^2)\cdots((2M-1)^2 - 4\nu^2)}{M!(8w)^M} + O(w^{-(M+1)})\right) \tag{7.7}$$

for $\nu = 0$ and 1. From the differential equation $wI_0''(w) + I_0'(w) - wI_0(w) = 0$ is obtained (for any $w \neq 0$), $I''(w) + I'(w)[2 + 1/w + I'(w)] + 1/w = 0$. Using this and its derivatives, together with (7.7), allows the construction of asymptotic series of the form $\frac{d^l I(w)}{dw^l} = \sum_{m=1}^M a_{m,l}w^{-m} + O(w^{-(M+1)})$ for various coefficients $a_{m,l}$. The existence of these series makes differentiation of (7.6) for $M = 1$ permissible (see Theorem 1.7.7 of Bleistein and Handelsman 1986), to give (7.6) for $M > 1$.

**Proof of Theorem 5.** From (2.8) and (2.2), for $2s \neq$ any $\gamma$

$$L_0'(s) = S_\chi(2s) + \sum_{\gamma<2s}\frac{d}{ds}I\left(\frac{2s}{\sqrt{1/4+\gamma^2}}\right) + \sum_{\gamma>2s}\frac{d}{ds}\log I_0\left(\frac{2s}{\sqrt{1/4+\gamma^2}}\right) \tag{7.8}$$

The second and third terms of (7.8) may be represented approximately by, respectively:

$$\frac{1}{2\pi}\int_{\gamma_0}^{2s}\log\left(\frac{qw}{2\pi}\right).\frac{2}{w}I'\left(\frac{2s}{w}\right)dw = \frac{1}{\pi}\int_{1}^{2s/\gamma_0}\left[\log\left(\frac{qs}{\pi}\right)-\log x\right]\frac{I'(x)}{x}dx$$

$$= \frac{1}{\pi}\left[\left[\log\left(\frac{qs}{\pi}\right)-\log x\right]\frac{I(x)}{x}\right]_1^{2s/\gamma_0} + \frac{1}{\pi}\int_{1}^{2s/\gamma}\frac{I(x)}{x^2}\left[\log\left(\frac{qs}{\pi}\right)+1-\log x\right]dx \quad (7.9)$$

$$\frac{1}{2\pi}\int_{2s}^{\infty}\log\left(\frac{qw}{2\pi}\right).\frac{2}{w}\left(\log I_0\left(\frac{2s}{w}\right)\right)'dw = \frac{1}{\pi}\int_{0}^{1}\left[\log\left(\frac{qs}{\pi}\right)-\log x\right]\frac{(\log I_0(x))'}{x}dx$$

$$= \frac{1}{\pi}\left[\left[\log\left(\frac{qs}{\pi}\right)-\log x\right]\frac{\log I_0(x)}{x}\right]_0^1 + \frac{1}{\pi}\int_{0}^{1}\left[\log\left(\frac{qs}{\pi}\right)+1-\log x\right]\frac{\log I_0(x)}{x^2}dx \quad (7.10)$$

where $\gamma_0 <$ all $\gamma$. In each case, the substitution $x = 2s/w$ has been made and an integration by parts carried out.

Both these representations introduce errors which are $O\left(\frac{\log s}{s}\right)$, as the following arguments show. (7.9) replaces $\frac{1}{\sqrt{1/4+\gamma^2}}$ by $1/\gamma$. Since $\frac{1}{w} - \frac{1}{\sqrt{1/4+w^2}} < \frac{1}{8w^3}$ the error has order

$$\frac{1}{16\pi}\int_{\gamma_0}^{2s}\frac{1}{w^3}\log\left(\frac{qw}{2\pi}\right)\frac{d}{dw}\left(\frac{2}{w}I'\left(\frac{2s}{w}\right)\right)dw = O\left(\frac{\log s}{s}\right)$$

integrating by parts and using (4.6) and (7.6). (7.9) also contains error from using (4.7) which is, integrating by parts twice:

$$\frac{1}{2\pi}\int_{\gamma_0}^{2s} N_R'(w).\frac{2}{w}I'\left(\frac{2s}{w}\right)dw = \frac{1}{\pi}\left[\frac{N_R(w)}{w}I'\left(\frac{2s}{w}\right)\right]_{\gamma_0}^{2s} - \frac{1}{\pi}\left[N_I(w)\frac{d}{dw}\left(\frac{1}{w}I'\left(\frac{2s}{w}\right)\right)\right]_{\gamma_0}^{2s}$$

$$+ \frac{1}{\pi}\int_{\gamma_0}^{2s} N_I(w)\frac{d^2}{dw^2}\left(\frac{1}{w}I'\left(\frac{2s}{w}\right)\right)dw = O\left(\frac{\log s}{s}\right)$$

assuming that $N_I(2s)$ as well as $N_R(2s)$ is $O(\log s)$, and again using (4.6) and (7.6). For large $s$ the error in replacing the upper limit by infinity in (7.9) is, using (7.6), of order

$$\frac{1}{\pi}\int_{2s/\gamma_0}^{\infty}\frac{1}{x^2}\left[\log s + \log\left(\frac{q}{\pi}\right) - \log x\right]dx = O\left(\frac{\log s}{s}\right)$$

The error in the representation (7.10) from replacing $\frac{1}{\sqrt{1/4+\gamma^2}}$ by $1/\gamma$ is, as for (7.9)

$$\frac{1}{8\pi}\int_{2s}^{\infty}\log\left(\frac{qw}{2\pi}\right).\frac{1}{w^4}\left(\log I_0\left(\frac{2s}{w}\right)\right)'dw = O\left(\frac{\log s}{s^3}\right)$$

using (5.1); and the error from using (4.7) is

$$\frac{1}{\pi}\int_{2s}^{\infty}\frac{N_R'(w)}{w}\left(\log I_0\left(\frac{2s}{w}\right)\right)' dw = \frac{-N_R(2s)I_1(1)}{2s\pi I_0(1)} - \frac{1}{\pi}\int_{2s}^{\infty}N_R(w)\frac{d}{dw}\left[\frac{1}{w}\left(\log I_0\left(\frac{2s}{w}\right)\right)'\right] dw$$

using (7.6); the integral is $O\left(\frac{\log s}{s}\right)$ using (5.1).

The result for $L_0'(s)$ in (2.19) follows on sorting out the various terms in (2.14), (7.9), and (7.10), noting that some terms involving $\log I_0(1)$ cancel. Validity when $2s =$ any $\gamma$ follows by continuity. The remaining results in (2.19) can be proved similarly. That for $L_0(s)$ could be obtained directly by integration of $L_0'(s)$, with the slightly less good error term $O((\log s)^2)$.

## 8. METHOD OF STEEPEST DESCENT TO THIRD ORDER

**Calculation of moments of transform distribution.** Theorem 5 can now be applied to a more accurate and rigorous steepest descent calculation than that leading to (7.2). Existence of the moment generating function means that all derivatives of $\bar{P}_0(s)$ exist, as for any $V$ do

$$\bar{p}_l(s, V) = \int_{-\infty}^{\infty}(v-V)^l e^{sv}P_0(v)dv = e^{sV}\int_{-\infty}^{\infty} v^l e^{sv}P_0(v+V)dv \tag{8.1}$$

$\bar{P}_0(s) = \bar{p}_0(s, V)$ and $\bar{P}_0'(s) = \int_{-\infty}^{\infty}ve^{sv}P(v)dv = V(s) + \bar{p}_1(s, V)$. So if as (7.1) $V(s) = L_0'(s)$, then $\bar{p}_1(s, V(s)) = 0$. $V(s)$ is the mean of the distribution $e^{sv}P_0(v)/\bar{P}_0(s)$, and

$$H_l(s) = \frac{\bar{p}_l(s, V(s))}{\bar{P}_0(s)} \tag{8.2}$$

are its moments. These can be calculated from $L_0(s)$ and its derivatives, since

$$\frac{d\bar{p}_l(s, V(s))}{ds} = \int_{-\infty}^{\infty}\left[v(v-V)^l - l(v-V)^{l-1}\frac{dV}{ds}\right]e^{sv}P_0(v)dv = \bar{p}_{l+1} + L_0'(s)\bar{p}_l - lL_0''(s)\bar{p}_{l-1}$$

$$H'_l(s) = -\frac{\bar{P}_0'(s)\bar{p}_l(s)}{(\bar{P}_0(s))^2} + H_{l+1}(s) + L_0'(s)H_l(s) - lL_0''(s)H_{l-1}(s) = H_{l+1}(s) - lL''(s)H_{l-1}(s)$$

Thus $H_{l+1}(s) = H'_l(s) + lL_0''(s)H_{l-1}(s)$, and so

$$H_2(s) = L_0''(s), \qquad H_3(s) = L_0'''(s),$$

$$H_4(s) = L_0''''(s) + 3(L_0''(s))^2, \qquad H_5(s) = L_0'''''(s) + 10L_0''(s)L_0'''(s) \tag{8.3}$$

**Relationship between $P$ and transform in terms of moments.** With the notation of (7.1)

$$\frac{H_l(s)}{e^{sV}} = \int_{-\infty}^{\infty} v^l e^{(sv+Q(V+v))}dv = P_0(V)\int_{-\infty}^{\infty} v^l \exp\left(v\Delta s - \frac{v^2}{2\eta^2} + \sum_{m=3}^{\infty}\frac{v^m}{m!}Q^{(m)}(V)\right)dv \tag{8.4}$$

Also define $\kappa(s) = \eta^3 Q'''(V(s))$ as before, and $\lambda(s)$, $\mu(s)$ and $\rho(s)$ by

$$\eta^4 Q''''(V) = \frac{Q''''(V)}{\left(-Q''(V)\right)^2} = \lambda \kappa^2; \quad \eta^5 Q'''''(V) = \frac{Q'''''(V)}{\left(-Q''(V)\right)^{5/2}} = \mu \kappa^3; \quad \eta \Delta s = \kappa \rho \qquad (8.5)$$

Assume for the present that $\lambda(s), \mu(s)$ and $\rho(s)$ are of order 1; results will confirm this. Retain terms to $m = 5$ in (8.4), work to third order in $\kappa$, and set $v = \eta w$. Expanding the exponential:

$$\bar{p}_l(s,V) \cong \eta^{l+1} e^{sV} P_0(V) \int_{-\infty}^{\infty} w^l e^{-w^2/2} \left[ 1 + \kappa w \left( \rho + \frac{w^2}{6} \right) + \frac{\lambda \kappa^2 w^4}{24} + \frac{\mu \kappa^3 w^5}{120} \right.$$
$$\left. + \frac{\kappa^2 w^2 (\rho + w^2/6)^2 + \lambda \kappa^3 w^5 (\rho + w^2/6)/12}{2} + \frac{\kappa^3 w^3 (\rho + w^2/6)^3}{6} \right] dw \quad (8.6)$$

So since $\bar{p}_1(s,V) = 0$

$$0 \cong \eta^2 e^{sV} \sqrt{2\pi} P_0(V) \kappa \left[ \rho + \frac{1}{2} + \kappa^2 \left( \frac{\mu}{8} + \frac{\rho^3}{2} + \frac{5\rho^2}{4} + \frac{35\rho}{24} + \frac{5\rho\lambda}{8} + \frac{35(1+\lambda)}{48} \right) \right],$$

which working to order $\kappa$ only gives $\rho \cong -1/2$ or $\Delta s \cong -\kappa/2\eta$ as in Section 7. To order $\kappa^3$ this value may be substituted in the terms above which are multiplied by $\kappa^2$, to give

$$\rho \cong -\frac{1}{2}\left( 1 + \kappa^2 \left( \frac{5\lambda}{6} + \frac{\mu}{4} + \frac{1}{2} \right) \right) \qquad (8.7)$$

Applying (8.7) to (8.6) for $l = 0, 2, 3, 4,$ and 5 gives, to order $\kappa^3$

$$\bar{P}_0(s) \cong \eta \, e^{sV} \sqrt{2\pi} P_0(V)[1 + \kappa^2(\lambda/8 + 1/12)]$$
$$\bar{p}_2(s,V) \cong \eta^3 e^{sV} \sqrt{2\pi} P_0(V)[1 + \kappa^2(5\lambda/8 + 7/12)]$$
$$\bar{p}_3(s,V) = \eta^4 e^{sV} \sqrt{2\pi} P_0(V)\kappa[1 + \kappa^2(25\lambda/8 + \mu/2 + 31/12)]$$
$$\bar{p}_4(s,V) \cong \eta^5 e^{sV} \sqrt{2\pi} P_0(V)[3 + 5\kappa^2(7\lambda/8 + 5/4)]$$
$$\bar{p}_5(s,V) = \eta^6 e^{sV} \sqrt{2\pi} P_0(V)\kappa[10 + \kappa^2(185\lambda/4 + 6\mu + 275/6)]$$

so that, dividing these results, to order $\kappa^3$

$$H_2(s) \cong \eta^2[1 + \kappa^2(\lambda+1)/2] \qquad H_3(s) \cong \eta^3 \kappa[1 + \kappa^2(3\lambda + \mu/2 + 5/2)]$$
$$H_4(s) \cong \eta^4[3 + 2\kappa^2(2\lambda + 3)] \qquad H_5(s) \cong \eta^5 \kappa[10 + \kappa^2(45(\lambda+1) + 6\mu)] \qquad (8.8)$$

**Solution for $P_0$: Proof of Theorem 6.** If are defined the ratios

$$L^{(3)}(s) = \frac{L_0'''(s)}{\left(L_0''(s)\right)^{3/2}} \; ; \quad L^{(4)}(s) = \frac{L_0''''(s)}{\left(L_0''(s)\right)^2} \; ; \quad L^{(5)}(s) = \frac{L_0'''''(s)}{\left(L_0''(s)\right)^{5/2}} \qquad (8.9)$$

then from (8.3) and (8.8), to order $\kappa^3$:

$$L^{(3)} \cong \kappa[1 + \kappa^2(9\lambda/4 + \mu/2 + 7/4)]; \; L^{(4)} \cong \kappa^2 (\lambda + 3); \; L^{(5)} \cong \kappa^3(10\lambda + \mu + 15)$$

Thus if for any $s$ $L^{(3)}, L^{(4)}$ and $L^{(5)}$ are calculated, then $\lambda, \mu$ and $\kappa$ can be determined by iterating

$$\lambda = L^{(4)}/\kappa^2 - 3, \; \mu = L^{(5)}/\kappa^3 - 10\lambda - 15, \; \kappa = L^{(3)}/[1 + \kappa^2(9\lambda/4 + \mu/2 + 7/4)] \qquad (8.10)$$

starting with the approximation $\kappa = L^{(3)}$. Then again from (8.3) and (8.8), to order $\kappa^3$:

$$\eta = \sqrt{L_0''(s)}\left(1 - \frac{\kappa^2(\lambda+1)}{4}\right) \quad ; \quad \Delta s \cong -\frac{\kappa}{2\eta}\left(1 + \kappa^2\left(\frac{5\lambda}{6} + \frac{\mu}{4} + \frac{1}{2}\right)\right) \tag{8.11},$$

and, for $V = L_0'(s)$, to order $\kappa^3$:

$$\log P_0(V) \cong L_0(s) - sV - \log\eta - \log(\sqrt{2\pi}) - \kappa^2(\lambda/8 + 1/12)$$

$$\cong L_0(s) - sV - \frac{1}{2}\log(2\pi L_0''(s)) + \kappa^2\left(\frac{\lambda}{8} + \frac{1}{6}\right) \tag{8.12}$$

Also for large $s$, from (2.19), similar results for higher derivatives of $L_0(s)$, (8.10) and (8.11):

$$\kappa \sim -\sqrt{\frac{\pi}{s\log s}}, \quad \lambda \sim -1, \quad \eta \sim \sqrt{\frac{\log s}{\pi s}}, \quad \Delta s \sim \frac{\pi}{2\log s}, \quad \mu \sim 1 \tag{8.13}$$

verifying earlier assumptions. The error in (8.12) is of order $\kappa^4$. That gives (2.23), and proves Theorem 6. The logarithmic error in (7.2) is of order $\kappa^2$, that is $O(1/(s\log s))$.

## 9. CALCULATIONS USING METHOD OF STEEPEST DESCENT

**Proof of Theorem 7.** This parallels exactly that of Theorem 4; $s$ replaces $i\omega$, giving (2.25). Throughout the range of convergence the series involved may be differentiated, giving estimates for the order of error in calculating derivatives of $L_u(s)$.

**Proof of Theorem 8.** This is given for a single series, but generalises. It is found convenient in multi-series cases to calculate $L_0(s)$ by summing the results of separate calculations for the various $\chi$. Write

$$L_u(s) = \Sigma(t) + \Sigma_E(t), \tag{9.1}$$

where $\quad \Sigma(t) = \displaystyle\sum_{k=1}^{\infty} \frac{(-1)^{k-1}(y + 1/(2k-1))t^{2k}}{k(2k-1)(y+1)j_1^2 T^{2(k-1)}}$

$$= \sum_{k=1}^{\infty} \frac{(-1)^{k-1}t^{2k}}{(y+1)j_1^2 T^{2(k-1)}}\left[\frac{2(y-1)}{2k-1} - \frac{y-1}{k} + \frac{2}{(2k-1)^2}\right]$$

which sums to the expression (2.26).

The proportionate error in (2.13) increases with $k$ no faster than $k$. From (2.17), the proportionate error in the approximation $c_k \cong 2^k/(kj_1^{2k})$ reduces faster than $(0.19)^k$. Therefore the terms of the series

$$\Sigma_E(t) = \sum_{k=1}^{\infty}(-1)^{k-1}\left(c_k r_k - \frac{y + 1/(2k-1)}{k(2k-1)(y+1)j_1^2 T^{2(k-1)}}\right)t^{2k} \tag{9.2}$$

eventually decrease more quickly than $(0.19\, t^2/T^2)^k$, completing the proof. A similar result could be applied to (2.16); however there is no point in doing so, since $\hat{P}_u(\omega)$ is negligible for values of $\omega$ such that (2.16) approaches divergence.

**Calculation of $L_0(s)$ for a particular $s$.** 30 terms of the power series $I_0(w) = \sum_{n=0}^{\infty} \frac{w^{2n}}{2^{2n}(n!)^2}$, $I_0'(w) = I_1(w) = \sum_{n=0}^{\infty} \frac{w^{2n+1}}{2^{2n+1}n!(n+1)!}$ give values to 1 part in $10^{14}$ for $|w| \leq 20$. 15 terms of the asymptotic expansions (7.7) give similar accuracy for $|w| \geq 20$. Thus $\log I_0(w)$ and its derivative $I_1(w)/I_0(w)$ can be calculated, and higher derivatives of $\log I_0(w)$ are expressible in terms of these using the differential equation $wI_0''(w) + I_0'(w) - wI_0(w) = 0$ and its derivatives. Thus the initial terms in (2.24) and its derivatives may be calculated. Values of $\Sigma(t)$ as given by (2.26) and expressions for its derivatives are easily calculated to high accuracy, the integral of inverse tangent by integrating the polynomial approximation (4.4.49) in Abramowitz and Stegun 1964. Convergence of (9.2) can be checked for each calculation: for all calculations here 8 terms were sufficient (with $c_k$ for $k > 5$ determined by 2 or 3 terms of (6.2)).

**Calculation of $P_0(v)$ for a particular $v$.** The value of $s$ such that $L_0'(s) = v$ can be obtained from an initial estimate $s_{(0)}$ by iteration. This is found to converge to the accuracy of computation in a few steps:

$$s_{(l+1)} = s_{(l)} + \left(v - L_0'(s_{(l)})\right)\frac{ds}{dv} = s_{(l)} + \frac{v - L_0'(s_{(l)})}{L_0''(s_{(l)})} \tag{9.3}$$

**Calculation of $E(v)$.** Expanding the integrand in $\frac{E(v)}{P_0(v)} = \frac{1}{P_0(v)}\int_0^{\infty} P_0(v + v')\,dv'$ as in (8.4) and (8.6):

$$\frac{E(v)}{P_0(v)} \cong E_0 + \frac{\kappa E_3}{6} + \kappa^2\left(\frac{\lambda E_4}{24} + \frac{E_6}{72}\right) + \kappa^3\left(\frac{\mu E_5}{120} + \frac{\lambda E_7}{144} + \frac{E_9}{1296}\right) \tag{9.4}$$

where $E_m = \eta \int_0^{\infty} w^m e^{-\eta(s-\Delta s)w - w^2/2}\,dw$.

For large $s$, $E_m < \eta \int_0^{\infty} w^m e^{-\eta s}\,dw = m!\,\eta/(\eta s)^{m+1}$. So from (8.13), the terms in (9.4) will be of order $1/s$, $\frac{1}{s^3(\log s)^2}$, $\frac{1}{s^4(\log s)^3}$, and $\frac{1}{s^5(\log s)^4}$ respectively. Given that the error in $\log P_0(v)$ as given by (2.23) is of order $\frac{1}{s^2(\log s)^2}$, all terms after the first should then be omitted, to give

$$\frac{E(v)}{P_0(v)} \cong E_0 \text{ where } E_0 = \frac{\eta\sqrt{2\pi}}{2}\exp\left(\frac{(\eta^2(s-\Delta s))^2}{2}\right)\left[1 - \text{erf}\left(\frac{\eta(s-\Delta s)}{\sqrt{2}}\right)\right] \tag{9.5}$$

When using (7.2) to first order in $\kappa$, $\eta \cong \sqrt{L_0''(s)}$, $\Delta s \cong -\kappa/2\eta \cong -L_0'''(s)/2(L_0''(s))^2$, and

$$\frac{E(v)}{P_0(v)} \cong \frac{1}{(s - \Delta s)}\left(1 - \frac{1}{\eta^2 s^2}\right) \cong \frac{1}{s}\left[1 - \frac{L_0'''(s)}{2s(L_0''(s))^2} - \frac{1}{s^2 L_0''(s)}\right] \tag{9.6}$$

Using (2.19), this gives

$$\frac{P_0(v)}{E(v)} \cong s\left[1 + \frac{\pi}{2s(\log s + A))}\left(1 - \frac{1}{\log s + A}\right)\left(1 + O\left(\frac{1}{s}\right)\right)\right] \tag{9.7}$$

whilst the logarithmic decrement of $P_0(v)$ is the slightly different

$$\frac{P_0'(v)}{P_0(v)} = \Delta s - s \cong \frac{L_0'''(s)}{2(L_0''(s))^2} - s$$

$$= -s\left[1 + \frac{\pi}{2s(\log s + A))}\left(1 + \frac{1}{\log s + A}\right)\left(1 + O\left(\frac{1}{s}\right)\right)\right] \quad (9.8)$$

For less large $s$, the full expression (9.4) for $E_0$ may be used, calculating $\eta$ and $\Delta s$ using (8.11), and obtaining other $E_m$ as necessary through integration by parts:

$$E_1 = \eta - \eta(s - \Delta s)E_0 \; ; \quad E_m = (m-1)E_{m-2} - \eta(s - \Delta s)E_{m-1} \quad (m \geq 2) \quad (9.9)$$

**Proof of Theorem 9.** This follows Monach 1980 and Lamzouri 2011, but working to $O(\log s)$ rather than $O(s)$. Where $q$ is prime and hence the values of $A$ are the same for all characters involved, the expressions for $L_0(s)$ and its derivatives are the sum of expressions of form (2.19) with $s$ replaced by $\alpha_\chi s$. So

$$L_0'(s) = \sum_\chi \alpha_\chi \Delta_\chi + \frac{\alpha_\Sigma}{2\pi}[(\log s)^2 + 2(A+Y)\log s + 2\pi X + 2AY + Z] + O\left(\frac{\log s}{s}\right)$$

$$L_0(s) = sL_0'(s) - \frac{s}{\pi}\{\log s + A + Y - 1\} + O(\log s)$$

and $v = L_0'(s)$ is a quadratic equation for $s$ whose positive solution is, with the notation of (2.28), $\log s = W - (A + Y) + O(1/s)$. Then (2.23) or (7.2) approximate to

$$\log P_0(v) = L_0(s) - sv + O(\log s) = -\frac{s}{\pi}\{\log s + A + Y - 1\} + O(\log s)$$

which gives (2.27), noting that $e^{-A} = \frac{\pi}{q}e^{-A_0}$, and that for large $s$, $\alpha_\Sigma(\log s)^2 < 2\pi v$. From (9.7), $\log E(v) - \log P(v) \sim \log s$, so (2.27) also gives $\log E(v)$ to the stated accuracy.

**Results of calculations.** Table 4 sets out results for $1 \leq v \leq 11$, for the prime count distribution; for square/non-square races modulo 4, 7 and 13; and for the race of residue 1 against 3 or 5 modulo 7, which involves 5 series. The most extreme results calculated correspond to the logarithmic density of $\pi(x) > \text{Li}(x) + 5\text{Li}(\sqrt{x})$ and to the results of the races going "against the trend" by a factor of 10 or more. Key points on the calculations are as follows:

a) For each $v$, calculations were done with two values of $u$, corresponding to a factor of about 2 difference in each $N(u)$ or $N_\chi(u)$ The aim was to keep the difference $\delta_u$ between values of $\log E(v)$ thereby obtained less than $5 \times 10^{-5}$. This corresponded to $t/T$ less than about 0.6, and the last term of (9.2) used ($k = 8$) less than $10^{-7}$. Once $v$ reached a value where $\delta_u > 5 \times 10^{-5}$, a new higher value of $N(u)$ was selected, the previous higher value becoming the lower. This reduced $\delta_u$ to $10^{-7}$ or less and confirmed accuracy up to each changeover. For the prime count distribution, calculations thus began by comparing 50 initial terms with 25, and concluded by comparing 16000 with 10000. The races required smaller numbers of initial terms – up to 3000 for modulo 4, up to 100 for each series for the inter – residue race.

TABLE 4. RESULTS FROM METHOD OF STEEPEST DESCENT

| $v$ | $\sqrt{2\pi v}$ | Prime count distribution | | | Mod 4 square/non-square | | | Mod 7 square/non-square | | | Mod 13 square/non-square | | | Mod 7, 1 vs. 3 or 5 | | |
|---|---|---|---|---|---|---|---|---|---|---|---|---|---|---|---|---|
| | | $\log E(v)$ | $\epsilon_2(\%)$ | $\epsilon_{ML}$ | $\log E(v)$ | $\epsilon_2(\%)$ | $\epsilon_{ML}$ | $\log E(v)$ | $\epsilon_2(\%)$ | $\epsilon_{ML}$ | $\log E(v)$ | $\epsilon_2(\%)$ | $\epsilon_{ML}$ | $\log E(v)$ | $\epsilon_2(\%)$ | $\epsilon_{ML}$ |
| 1.0 | 2.51 | -15.1511 | -0.219 | 1.83 | -5.504 | -1.090 | 0.23 | | | | | | | | | |
| 1.5 | 3.07 | -35.9677 | -0.091 | 1.91 | -11.483 | -0.482 | 0.08 | -7.422 | -1.149 | 0.18 | | | | | | |
| 2.0 | 3.54 | -72.1428 | -0.045 | 1.97 | -21.5027 | -0.221 | -0.06 | -13.327 | -0.458 | 0.04 | -9.070 | -0.726 | -0.27 | | | |
| 2.5 | 3.96 | -130.1514 | -0.025 | 2.00 | -37.2592 | -0.115 | -0.19 | -22.5157 | -0.231 | -0.09 | -14.834 | -0.417 | -0.45 | | | |
| 3.0 | 4.34 | -218.2738 | -0.015 | 2.03 | -60.9113 | -0.066 | -0.31 | -36.2276 | -0.129 | -0.21 | -23.294 | -0.236 | -0.62 | | | |
| 3.5 | 4.69 | -346.9825 | -0.010 | 2.05 | -95.1814 | -0.040 | -0.43 | -56.0264 | -0.077 | -0.32 | -35.3599 | -0.138 | -0.77 | -11.301 | -0.244 | 1.64 |
| 4.0 | 5.01 | -529.3893 | -0.006 | 2.06 | -143.4698 | -0.026 | -0.54 | -83.8639 | -0.049 | -0.43 | -52.1624 | -0.087 | -0.92 | -14.514 | -0.239 | 1.59 |
| 4.5 | 5.32 | -781.7636 | -0.004 | 2.08 | -209.9886 | -0.018 | -0.64 | -122.1566 | -0.032 | -0.54 | -75.0992 | -0.057 | -1.07 | -18.3185 | -0.205 | 1.54 |
| 5.0 | 5.60 | -1124.1332 | -0.003 | 2.09 | -299.9155 | -0.012 | -0.75 | -173.8739 | -0.022 | -0.64 | -105.8832 | -0.039 | -1.20 | -22.7800 | -0.164 | 1.50 |
| 5.5 | 5.88 | -1580.9784 | -0.002 | 2.09 | -419.5716 | -0.009 | -0.84 | -242.6404 | -0.016 | -0.74 | -146.6006 | -0.027 | -1.33 | -27.9695 | -0.128 | 1.46 |
| 6.0 | 6.14 | -2182.0308 | -0.002 | 2.10 | -576.6257 | -0.006 | -0.94 | -332.8527 | -0.011 | -0.83 | -199.7766 | -0.020 | -1.46 | -33.9624 | -0.101 | 1.42 |
| 6.5 | 6.39 | -2963.1926 | -0.001 | 2.11 | -780.3300 | -0.005 | -1.03 | -449.8143 | -0.008 | -0.92 | -268.4518 | -0.015 | -1.58 | -40.8393 | -0.082 | 1.38 |
| 7.0 | 6.63 | -3967.5898 | -0.001 | 2.11 | -1041.7884 | -0.004 | -1.12 | -599.8896 | -0.006 | -1.01 | -356.2690 | -0.011 | -1.70 | -48.6867 | -0.068 | 1.34 |
| 7.5 | 6.86 | -5246.7766 | -0.001 | 2.11 | -1374.2653 | -0.003 | -1.20 | -790.6805 | -0.005 | -1.10 | -467.5724 | -0.008 | -1.81 | -57.5969 | -0.058 | 1.31 |
| 8.0 | 7.09 | -6862.1110 | -0.001 | 2.12 | -1793.5359 | -0.002 | -1.28 | -1031.2272 | -0.004 | -1.18 | -607.5209 | -0.006 | -1.92 | -67.6691 | -0.050 | 1.27 |
| 8.5 | 7.31 | -8886.3200 | 0.000 | 2.12 | -2318.2868 | -0.002 | -1.36 | -1332.2374 | -0.003 | -1.26 | -782.2170 | -0.005 | -2.03 | -79.0088 | -0.043 | 1.23 |
| 9.0 | 7.52 | -11405.2800 | 0.000 | 2.12 | -2970.5695 | -0.001 | -1.44 | -1706.3463 | -0.002 | -1.34 | -998.8531 | -0.004 | -2.13 | -91.7288 | -0.037 | 1.20 |
| 9.5 | 7.73 | -14520.0358 | 0.000 | 2.12 | -3776.3165 | -0.001 | -1.52 | -2168.4118 | -0.002 | -1.42 | -1265.8767 | -0.003 | -2.24 | -105.9495 | -0.032 | 1.17 |
| 10.0 | 7.93 | -18349.0874 | 0.000 | 2.12 | -4765.9238 | -0.001 | -1.60 | -2735.8482 | -0.001 | -1.49 | -1593.1788 | -0.002 | -2.34 | -121.7989 | -0.028 | 1.13 |
| 10.5 | 8.12 | -23030.9749 | 0.000 | 2.13 | -5974.9105 | -0.001 | -1.67 | -3429.0042 | -0.001 | -1.56 | -1992.3049 | -0.002 | -2.43 | -139.4131 | -0.024 | 1.10 |
| 11.0 | 8.31 | -28727.1968 | 0.000 | 2.13 | -7444.6626 | -0.001 | -1.74 | -4271.5880 | -0.001 | -1.64 | -2476.6945 | -0.002 | -2.53 | -158.9371 | -0.021 | 1.07 |

NOTES. $\log E(v)$ is calculated from (2.23) with (9.5) and if necessary (9.4) and (9.8). $\epsilon_2(\%) = 100 \times$ difference in $\log E(v)$ as thus calculated and the simpler calculation using (7.2) and (9.7) (thus = proportionate difference between the two calculations of $E(v)$). $\epsilon_{ML}$ = absolute difference between $\log E(v)$ as calculated from (2.23) with (8.10), and as estimated using (2.27) and (2.28).

b) The derivative of (2.24) represents $v$ as the sum of contributions from the initial terms and the remainder. Typically the remainder contributes about 30% of smaller $v$, falling to 20% of larger $v$; and the total contribution of initial terms to $v$ is 70 – 90% of their maximum possible total contribution, $S(u)$. Thus typically a value of $N(u)$ is needed such that $S(u) \approx 0.8v$. The smaller numbers of initial terms needed for races than for the prime count distribution reflects the larger values of $S(u)$ for a particular $u$.

c) The iterative process to solve (8.10) converges rapidly, to (for larger $v$) values of $\lambda$ between -0.8 and -0.9, and values of $\mu$ between 0.6 and 0.7, tending towards the limiting values -1 and 1 respectively.

d) The accurate calculation of the remainder $\sum_{k=1}^{\infty}(-1)^k c_k r_k t^{2k}$ is very necessary for large $v$. Even omitting the first term (corresponding to a normal distribution) it can contribute several hundred to $\log E(v)$.

e) The correction (shown as a percentage in Table 4) from using (2.23) rather than (7.2) diminishes as $v$ increases, and for the prime count distribution is well under 0.01% for $v > 5$. Because of the negative value of $\lambda$ it is generally less in magnitude than $0.2(s \log s)^{-1}$.

f) As (2.23) has error $O((s \log s)^{-2})$, its accuracy can be expected to improve with increasing $v$ and $s$. The races with larger $\sigma_0$ (see Table 2) have less extreme results for a specified $v$, and hence a larger $v$ is necessary to attain high accuracy. This is discussed further in Section 12 where comparisons with other methods are considered, in the ranges of $v$ where these are capable of providing accurate answers. Results in Table 4 are quoted to 3 rather than 4 decimal places where these comparisons indicate only that accuracy is justified, and are omitted if fewer than 3 places are justified.

g) (2.27) gives results differing from calculated values by about 2. This is of order $\log s$, though corresponding to a factor of order 10 difference from the calculated $E(v)$. The errors are systematic, suggesting that further improvements in accuracy may be possible. This has not been pursued, since (2.27) is sufficient to explain major aspects of behaviour of $E(v)$; also, uncertainties in the estimates of $\Delta_\chi$ made from the data on zeros available could change the quoted differences from calculated values by up to 0.2 – 0.3 at $v = 11$. The simpler formula of Monach and Lamzouri, that is $\log P_0(v) \sim -\frac{\sqrt{2\pi v}}{q}\exp(\sqrt{2\pi v} - A_0)$ has error of order $s$ rather than $\log s$, and underestimates $E(11)$ by a factor of order $10^{3630}$.

All this gives high confidence, though not rigorous proof, that the calculated values of $\log E(v)$ are accurate as stated, and that $E(v)$ is therefore accurate to 0.01% for larger $v$.

## 10. EXPLICIT FORM OF $P_u(v)$ FOR SMALL $v$: APPLICATION TO NON-EXTREME RACES

An analogue of the steepest descent method may be applied to invert the Laplace transform, thereby obtaining an explicit expression for the deviation of $P_u(v)$ from normal.

**Proof of Theorem 10.** For $< T$, setting $Y = v/\sigma_u$

$$\exp\left(\sum_{k=1}^{\infty}(-1)^{k-1} c_k R_k t^{2k}\right) = \bar{P}_u\left(\frac{Y}{\sigma_u}\right) = \int_{-\infty}^{\infty} P_u(v') \exp\left(\frac{v'Y}{\sigma_u}\right) dv'$$

$$= \int_{-\infty}^{\infty} P_u(v') \exp\left(\frac{v'^2}{2\sigma_u^2} - \frac{1}{2}\left(\frac{v'}{\sigma_u} - \Upsilon\right)^2 + \frac{\Upsilon^2}{2}\right) dv' = \frac{e^{\Upsilon^2/2}}{\sqrt{2\pi}} \int_{-\infty}^{\infty} G_u(\Upsilon + \Upsilon') e^{-\Upsilon'^2/2} d\Upsilon'$$

where $\Upsilon' = v'/\sigma_u - \Upsilon$. Thus

$$\exp\left(\sum_{k=2}^{\infty}(-1)^{k-1} c_k R_k \Upsilon^{2k}\right) = \frac{1}{\sqrt{2\pi}} \int_{-\infty}^{\infty} G_u(\Upsilon + \Upsilon') e^{-\Upsilon'^2/2} d\Upsilon'$$

$$= \frac{1}{\sqrt{2\pi}} \int_{-\infty}^{\infty} e^{-\Upsilon'^2/2} \sum_{l=0}^{\infty} \frac{\Upsilon'^l}{l!} \left(\frac{d}{d\Upsilon}\right)^l G_u(\Upsilon) d\Upsilon' = \sum_{m=0}^{\infty} \frac{(2m)!}{2^m m!} \cdot \frac{1}{(2m)!} \left(\frac{d}{d\Upsilon}\right)^{2m} G_u(\Upsilon) = \exp\left(\frac{1}{2}\frac{d^2}{d\Upsilon^2}\right) G_u(\Upsilon)$$

Inverting this operational series gives a result which is easily seen to be divergent, but indicates that for suitable $M$ and $K$ the asymptotic representation (2.30) should approximate $G_u(\Upsilon)$.

**Approximation to first order in $R_2$ for the deviation from normal.** For small $R_2$ and $\Upsilon/T$ just the $k = 2$ term can be retained in (2.30) and the $m$ series terminated at $m = 2$, giving:

$$P_u(v) \cong \frac{1}{\sigma_u\sqrt{2\pi}}\left(1 - \frac{3R_2}{16} + \frac{3R_2 v^2}{8\sigma_u^2}\right) \exp\left(-\frac{v^2}{2\sigma_u^2} - \frac{R_2 v^4}{16\sigma_u^4}\right) \qquad (10.1)$$

When $R_2 \Upsilon^4/16 \ll 1$, a further approximation is

$$G_u(\Upsilon) \cong (1 - 3R_2/16 + 3R_2\Upsilon^2/8 - R_2\Upsilon^4/16) \qquad (10.2)$$

$$P_u(v) \cong \frac{1}{\sigma_u\sqrt{2\pi}}\left(1 - \frac{3R_2}{16} + \frac{3R_2 v^2}{8\sigma_u^2} - \frac{R_2 v^4}{16\sigma_u^4}\right) \exp\left(-\frac{v^2}{2\sigma_u^2}\right) \qquad (10.3)$$

which shows the deviation of $P_u$ from normal form. The density is reduced by proportion $3R_2/16$ at $\Upsilon = 0$. The equals normal density at $\Upsilon = \sqrt{3 - \sqrt{6}} \cong 0.742$, is increased in proportion $3R_2/8$ at $\Upsilon = \sqrt{3} \cong 1.732$, equals normal density again at $\Upsilon = \sqrt{3 + \sqrt{6}} \cong 2.334$, then falls away steadily.

This is the type of behaviour observed for the prime count distribution (see Figure 1), for which $P_0(v)$ might be validly represented by (2.30) when $\Upsilon < j_1 \sigma \sqrt{1/4 + \gamma_1^2}/2 \cong 3.65$ ($v \cong 0.78$) (10.1) and (10.3) could be expected to break down for some smaller $v$. In fact, (10.1) approximates the exact values of $P_0(v)$, as calculated using (6.9), to within 1.5% right up to $v = 0.78$, and the approximation continues fair to around $v = 1$ ($\Upsilon \cong 4.65$). Numerical integration of (10.1) gives $E(1) \cong 2.45 \times 10^{-7}$, within 7% of the accurate result.

$G_0(v)$ represents $P_0(v)$ well for $0 < v < 1$ in the case of an inter-residue race when $q$ is large (making $\sigma_0$ large, and $\beta$ small in accordance with Theorem 2.) It is then valid to make the further approximation $P_0(v) \cong \frac{1}{\sigma_0\sqrt{2\pi}}\left(1 - \frac{3\beta}{16} - \frac{v^2}{2\sigma_0^2}\right)$, which gives (2.31). For $q = 43$, results calculated using (2.31) are given with Figure 4. These have a maximum error of about 0.0001, whilst assuming the distribution is normal has a maximum error of about 0.001.

## 11. CALCULATIONS BY CONVOLUTION

In calculation it is found that the terms of the asymptotic series alternate in sign, so that the apparent error is less than the last term included. Assuming this is generally true, for the prime count distribution

$$G_u(\Upsilon) \cong \sum_{m=0}^{5} \frac{1}{2^m m!} \left(-\frac{d^2}{d\Upsilon^2}\right)^m \exp\left(\frac{\Upsilon^2}{2} - \Sigma(\Upsilon) - \Sigma_E(\Upsilon)\right) \qquad (11.1)$$

with $\Sigma(\Upsilon)$ and $\Sigma_E(\Upsilon)$ as in (2.26) and (9.2) and $K = 5$, is found to give values of $G_u$ accurate to 0.01% out to $\Upsilon = 20$ with $N(u) = 250$ (T $\cong$ 33.98) and to $\Upsilon = 30$ with $N(u) = 800$ (T $\cong$ 58.32). These are consistent with an alternative derivation of (2.30), made by shifting through $i\Upsilon$ the contour of integration in the inversion $P_0(v) = \frac{1}{2\pi}\int_{-\infty}^{\infty} \hat{P}_0(\omega)e^{i\omega v}d\omega$, with $\hat{P}_0(\omega)$ when significant represented by (6.3), which suggests that (11.1) is valid for $|\Upsilon| < 0.6T$. For smaller $\Upsilon$, accuracy is much higher. Derivatives of $\Sigma(\Upsilon)$ and $\Sigma_E(\Upsilon)$ up to the fourth can be calculated as in Section 9, and the first 3 terms of (11.1) thus worked out: the remaining terms, which have smaller effect, can be obtained sufficiently accurately by numerical differentiation.

Define average values $\tilde{f}_\gamma(v)$ at values of $v$ separated by a step length $\Delta v$ as

$$\tilde{f}_\gamma(v) = \frac{1}{\pi \Delta v}\left[\arccos\left(2(v - \Delta v/2)\sqrt{1/4 + \gamma^2}\right) - \arccos\left(2(v + \Delta v/2)\sqrt{1/4 + \gamma^2}\right)\right] \qquad (11.2)$$

with suitable adjustments near end points $v = \pm 2/\sqrt{1/4 + \gamma^2}$. Then an approximation $\tilde{F}_u(v)$ to $F_u(v)$ can be built up by repeated convolution

$$\tilde{F}_u(v) = \Delta v \sum_{v'} \tilde{F}_{u'}(v - v')\tilde{f}_\gamma(v') \qquad (11.3)$$

where there is a single $\gamma$ with $u' < \gamma < u$. An approximation to $P_0(v)$ is then

$$P_0(v) = \Delta v \sum_{v'} \tilde{F}_u(v - v')G_u(v')e^{-v'^2/2\sigma_u^2}. \qquad (11.4)$$

A small $\Delta v$ is needed, for $\tilde{f}_\gamma$ to be a good representation near end points. Inaccuracies there translate into growing inaccuracy of $\tilde{F}_u(v)$ as $v$ nears $S(u)$. However, a lower limit on $\Delta v$ is enforced both by the computation time (which behaves as $\Delta v^{-2}$) and by the differencing error introduced in (11.2). On a desktop computer in single precision, the limit appears to be about 0.00001. Such a small $\Delta v$ is needed only in working out the distributions of "blocks" of typically 50 terms. These distributions are smooth, and can be combined accurately by numerical convolution with a larger $\Delta v$ ($\Delta v = 0.0001$ was used); the same is true of (11.4) and of a numerical integration of $P_0(v)$ to estimate $E(v)$. (11.4) needs to be summed only over the range which contributes significantly, effectively that where logarithmic derivatives of $\tilde{F}_u$ and $G_u$ are similar in magnitude (as is the basis of the method of steepest descent.)

For the prime count distribution, calculations were done with several $\Delta v$. As for the method of steepest descent, $u$ must be sufficient that $S(u)$ is not much less than $v$. The method is fairly easy to set up, but is heavy on computing time: with $N(u) = 800$, to obtain values of $P_0(v)$ and $E(v)$ to $v = 5$ in steps of 0.0001 for the prime count distribution required about 3 hours.

With $N(u) = 250, 400,$ and $600$, agreement with the 800 term result to within 0.01% was obtained for $v \leq 3.8, 4.25$ and $4.65$ respectively.

Calculations were also done for some races. As with other methods, the need to bring together more series was balanced by a reduced need for numbers of initial terms in each. With the results sought being less extreme, at least 6 decimal places can be obtained with $N(u) = 10$ at the most, and $\Delta v = 0.0005$  However, the Rubinstein-Sarnak method as developed in Section 6 of this paper is generally more accurate as well as needing less computing time, because the representation of the transform of the remainder is more accurate than (11.1).

One particular application is to multi-way races with real characters as considered in Feuerverger and Martin 2000, Martin 2000. When $q = 8$, labelling the densities and probabilities of excess of the distributions as $P_{-8}(v)$ and $E_{-8}(v)$ for $X_{-8}$ associated with Dirichlet character $\chi_{-8}$, similarly for $\chi_{-4}$ and $\chi_8$, then the order of residues will be 7, 1, 5 when $X_{-8} + X_8 > 2$ and $X_{-8} + X_{-4} < 2$. For any value $w$ of $X_{-8}$, the probabilities of these events will be $E_8(2-w)$ and $1 - E_{-4}(2-w) = E_{-4}(w-2)$ respectively, so the race result is

$$\int_{-\infty}^{\infty} P_{-8}(w) E_{-4}(w-2) E_8(2-w) dw \tag{11.5}$$

Similarly, the result of the race 3, 5, 1, 7 ( $X_{-8} - X_{-4} > 0$, $X_{-8} + X_8 < 2$, $X_{-8} + X_{-4} > 2$ ) is

$$\int_{1}^{\infty} P_{-8}(w)(E_{-4}(2-w) - E_{-4}(w)) E_8(w-2) dw \tag{11.6}$$

and the result of the square/non-square race (1 leads each of 3, 5 and 7) is (by addition)

$$\int_{1}^{\infty} [P_{-8}(w-2)E_{-4}(w)E_8(w) + P_{-4}(w-2)E_8(w)E_{-8}(w) + P_8(w-2)E_{-8}(w)E_{-4}(w)] dw$$
$$+ \ E_{-8}(1)E_{-4}(1)E_8(1) \tag{11.7}$$

These and similar results may be evaluated directly by convolution, with only one numerical integration, rather than needing a multidimensional Poisson sum as (2.57) of Feuerverger and Martin 2000. (The method may be shown equivalent to that sum by expressing the distributions as inverse transforms and using $\delta(v) = \frac{1}{2\pi}\int_{-\infty}^{\infty} e^{iwv} dw$.) With $N(u) = 10$ for each character, and $\Delta v = 0.0005$, accuracy of $10^{-7}$ is obtained. The method may also be applied to some races with $q = 24$, the $P$s and $E$s then representing the combined distributions associated with 3 pairs of the 7 characters involved. An analogue of (11.7) estimates the result of the race 1 leading each of 5, 7 and 11 as about $2.86 \times 10^{-8}$. The full square/non-square result (1 leading these and also 17, 19 and 23) must be less than this, so this race is more extreme than the pi vs. Li race.

## 12. COMPARISON OF RESULTS

**Comparison of the three methods.** As shown in Table 5, for the prime count distribution within $1 \leq v \leq 1.2$ there is agreement to within 0.001% amongst the three methods considered here of applying a better approximation for the remainder: Rubinstein-Sarnak, steepest descent (to third order) and convolution. Similar comparisons arise for other distributions, though at larger values of $v$, reflecting larger $\sigma_0$ and smaller values of $s$ corresponding to a particular $v$.

The limitation of about $10^{-16}$ in absolute accuracy in the Rubinstein – Sarnak calculation (as applied here) restricts its application to $|v/\sigma_0|$ less than about 6. Thus for the prime count distribution loss of accuracy is noticeable at $v = 1.3$ and is significant at $v = 1.4$.

Results of calculations by convolution gradually increase above those of the method of steepest descent, the difference reaching 0.2% at $v = 5$. Bearing in mind that the accuracy of the steepest descent method increases with $v$, the likely cause of this difference is (11.3) overestimating $\tilde{F}_u(v)$ on account of the "cup shaped" nature of $f(x) = \frac{1}{\pi\sqrt{1-x^2}}$. Convolutions with larger $\Delta v$ are less accurate, which bears this out.

For small $v$ the steepest descent method loses accuracy, through the impact of terms higher than third order in the parameter $\kappa \sim -\sqrt{\frac{\pi}{s\log s}}$. Where comparisons with accurate calculations by other methods are available, the error in this method is typically of magnitude around $0.1\kappa^4 \approx (s\log s)^{-2}$. In assessing accuracy of this method for larger $v$, it has been assumed that the error continues at such magnitude. The fact that the correction $\kappa^2\left(\frac{\lambda}{8} + \frac{1}{6}\right)$ from (7.2) to (2.23) is typically of magnitude around $0.1\kappa^2$ gives some support to this assumption.

TABLE 5. PRIME COUNT DISTRIBUTION: % DIFFERENCES IN CALCULATED $E(v)$

| $v$ | 0.8 | 0.9 | 1.0 | 1.1 | 1.2 | 1.3 | 1.4 | 2.0 | 3.0 | 4.0 | 5.0 |
|---|---|---|---|---|---|---|---|---|---|---|---|
| Convol/s-d | 0.010 | 0.003 | 0.000 | -0.001 | -0.001 | 0.000 | 0.000 | 0.003 | 0.015 | 0.053 | 0.195 |
| R-S/s-d | 0.010 | 0.002 | 0.000 | -0.001 | -0.001 | -0.008 | -0.544 | | | | |
| Convol/R-S | 0.000 | 0.000 | 0.000 | 0.001 | 0.001 | 0.008 | 0.544 | | | | |

The Rubinstein-Sarnak method is the best way to obtain results which are much larger than the absolute accuracy of computation. The method of steepest descent obtains accurate results for extreme deviations, but requires a separate, iterative calculation for each value of $v$ selected. Convolution can deliver a full distribution with fair accuracy over a range of values up to something fairly but not very extreme.

**Comparison with observed distributions.** As already noted, calculated $P_0(v)$ show the type of behaviour observed in Figure 1 for observed actual values of $D(x)$ as given by (3.2). However, the distribution of the actual values deviates above and below normal to a greater extent than calculation suggests. The comparison of moment ratios for the prime count distribution in Table 6 confirms this. Whilst the actual values have variance close to its limiting value, higher moments are significantly lower. The difference is not dependent on the exact start and end points of the range, or the exact logarithmic interval used. So it is not a consequence of the "end effects" which are ignored in deriving (1.1), or of the omission of extreme values of $D$ in sampling. Nor is it on a clear reducing trend with increasing $v$, as would be consistent with it arising from terms neglected in (1.1). In fact, the tables provided by Kulsha to very fine resolution around extreme values suggest that $D < -0.8$ for a logarithmic range of about 0.0003 within the above range of $x$. Assuming a symmetric distribution, this is consistent with $E(0.8) \cong 4.4 \times 10^{-6}$, compared to $4.72 \times 10^{-5}$ calculated from (6.8). The limiting distribution is established only over a wider range of $x$ than that for which $\pi(x)$ and $D$ can be calculated exactly.

TABLE 6. MOMENT RATIOS

| $k$ | 2 | 4 | 6 | 8 |
|---|---|---|---|---|
| $M_{2k}/M_{2k}^0$ – for actual $D$, $x = 10^5$ to $10^{10}$ | 1.0077 | 0.9182 | 0.7638 | 0.5822 |
| $M_{2k}/M_{2k}^0$ – for actual $D$, $x = 10^{10}$ to $10^{15}$ | 0.9749 | 0.9260 | 0.8385 | 0.7100 |
| $M_{2k}/M_{2k}^0$ – for actual $D$, $x = 10^{15}$ to $10^{20}$ | 0.9977 | 0.9217 | 0.8208 | 0.7076 |
| $M_{2k}/M_{2k}^0$ - from (4.5) | 1.0000 | 0.9652 | 0.9034 | 0.8229 |

Table 7 examines available data from Bays and Hudson 1999 (B) and Demichel 2005 on what may be the first excursions of $D$ above values of $v$ up to 1.5. The excursion occupies a logarithmic interval $\Delta z = \log(x_2/x_1)$, where $x_1$ and $x_2$ are its beginning and end points. Values of these are obtained by measurement from the graphs these authors provide. This process is very approximate, but $\Delta z/z$ should give an indication of the logarithmic frequency of excess, subject to a large sampling error. Its values agree with those of $E(v)$ as calculated above to within a reasonable factor of about 3. No "wildly improbable" excursion occurs, nor does an "expected" excursion not appear. This suggests that if analysis were carried out over very long ranges similar to that done by Bays and Hudson for $v = 1$, similar agreement with the values of $E(v)$ could be obtained.

TABLE 7. POSSIBLE FIRST EXCURSIONS (APPROXIMATE VALUES OF $\Delta z/z$)

| $v$ | Location $x$ | $z = \log x$ | $\Delta z$ | $\Delta z/z$ | $E(v)$ |
|---|---|---|---|---|---|
| 1.2 | $1.6 \times 10^{9608}$ | $2.2 \times 10^4$ | 0.00001 | $5 \times 10^{-10}$ | $2.833 \times 10^{-10}$ |
| 1.25 | $1.3 \times 10^{651157}$ | $1.5 \times 10^6$ | 0.0004 | $3 \times 10^{-11}$ | $3.839 \times 10^{-11}$ |
| 1.3 | $6.6 \times 10^{30802655}$ | $7.1 \times 10^7$ | 0.0002 | $3 \times 10^{-12}$ | $4.587 \times 10^{-12}$ |
| 1.5 | $5.5 \times 10^{1625185852}$ | $3.7 \times 10^{10}$ | 0.00003 | $8 \times 10^{-16}$ | $2.396 \times 10^{-16}$ |

## 13. CONCLUSIONS AND SUGGESTIONS FOR FURTHER WORK

The following general conclusions on the behaviour of $P_0(v)$ and $E(v)$ represent the expression in "$v$ space" (Theorems 3, 6, 9 and 10) of results on moments (Lemma 1), on the Fourier Transform (Theorem 4) and on the Laplace Transform (Theorems 5 and 7). They rest on the behaviour of parameters of the distribution as given by Theorems 1 and 2.

i. For values of $v$ up to a few standard deviations $\sigma_0$, $P_0(v)$ differs from normal by being lower for small $v$, higher when $v$ has values around $\sigma_0$ to $2\sigma_0$, then lower again. The key determinant of the difference is $\beta$. Inter-residue race distributions involve more and more separate series and have small values of $\beta$, becoming smaller with increasing $q$; they tend to normal with increasing $q$ at any given $v/\sigma_0$. Square/non-square race distributions have values of $\beta$ which tend to increase with $q$, leading eventually to a double humped distribution, far from normal. This arises because of the appearance of smaller and smaller initial zeros $\gamma_1$ (see Bays *et al* 2001.) Values of $E(1)$ can be calculated using small numbers of zeros explicitly, the rest being included in a collective remainder term of form (2.16), or (2.30) with (2.29).

ii. For large values of $v/\sigma_0$, the distribution behaves as (2.27). Its main driving factor is the long term behaviour of zeros as given by (4.7). Behaviour is influenced by $q$ and by the $\Delta_\chi$ , which are effectively corrections to $v$ determined largely by the initial zeros.

Thus the initial zeros affect both the deviation from normal for smaller $v/\sigma_0$ and the adjustments to (2.27) for larger $v/\sigma_0$, though in different ways. Lamzouri (2012) shows that for inter-residue races with large $q$, and hence small $\beta$ , this type of behaviour dominates once $v/\sigma_0 \gg \log q$ (his Theorem 4, noting that $\sigma_0 \sim \sqrt{\phi(q) \log q}$).

**Possibilities for further work.** These include

a) Application of the Rubinstein-Sarnak method with modified remainder, (6.9) with (2.16), to wider categories of prime number race, including races involving many series and multi-way races, for which the use of small $N_\chi$ should be particularly advantageous in reducing computing time.
b) Further computations and analytical work applying the method of steepest descent (2.23) to wider categories of race, including ascertaining whether the error in (2.27) can be further reduced.
c) More rigorous demonstration of some of the results. Presently it is the close agreement between separate computations, and between computation and theory, that primarily supports accuracy. Orders of magnitude for error, but not explicit limits, are given.

## 14. ACKNOWLEDGEMENTS

I am grateful to Daniel Fiorilli for helpful suggestions, including concerning the use of the method of steepest descent, to Andrey Kulsha for guidance on the use of his tables, and to Peter Dixon, Greg Martin and Trevor Wooley for encouragement.

## 15. REFERENCES

M.Abramowitz and I.A.Stegun, *Handbook of Mathematical Functions*, US National Bureau of Standards, 1964.

C.Bays and R.H.Hudson, "*Zeros of Dirichlet L–functions and irregularities in the distribution of primes*" Mathematics of Computation, **69** (1999), 861 – 866

C.Bays and R.H.Hudson, "*A new bound for the smallest $x$ with $\pi(x) > Li(x)$*", Mathematics of Computation, **69** (1999), 1285 – 1296

C. Bays, K.Ford, R.H.Hudson and M. Rubinstein, "*Zeros of Dirichlet L – functions near the Real Axis and Chebyshev's Bias*", Journal of Number Theory, **87** (2001), 54 – 76

N.Bleistein and R.A.Handelsman, *Asymptotic Expansions of Integrals,* Dover Books 1986

A.Feuerverger and G.Martin, "*Biases in the Shanks-Renyi prime number race*", Experimental Mathematics, **9** (2000), 535 – 570.

P. Demichel, http://sites.google.com/site/dmlpat2/Li_crossover_pi.pdf..

D. Fiorilli 2012, "*Highly Biased Prime Number Races*", arXiv 1210.6946

D. Fiorilli and G.Martin 2012, "*Inequities in the Shanks-Renyi Prime Number Race: an asymptotic formula for the densities*" www.math.ubc.ca/~gerg

A.E.Ingham, *The Distribution of Prime Numbers*, Cambridge 1932

A.V. Kulsha, "*The Fluctuations of the Prime-counting Function pi(x)*",
www.primefan.ru/stuff/primes/table.html

Y. Lamzouri, "*Large deviations of the Limiting distribution in the Shanks-Renyi prime number race*", Math. Proc. Cambs. Phil. Soc, **153**, 01 (July 2012), 147-166

P.Leboeuf, "*Prime Correlations and Fluctuations*", Annales Henri Poincare, **4** (2003), S727 – S752

R.Sherman Lehman, "*On the difference $\pi(x) - Li(x)$*", Acta Arithmetica, **XI** (1966), 397 - 410

J.E.Littlewood, "*Sur la distribution des nombres premiers*", C.R.Acad.Sci Paris, **158** (1914), 1869-1872

G.Martin "*Asymmetries in the Shanks – Renyi Prime Number Race*", Proceedings of the Millennial Conference on Number Theory, University of Illinois, **II**, (2000),403 – 415.

W.R.Monach, "*Numerical investigation of several problems in number theory*", University of Michigan Ph. D, Thesis (1980, unpublished).

H.L.Montgomery, "*The zeta function and prime numbers*", Proceedings of the Queen's Number Theory Conference, 1979, 1-31

H.L.Montgomery and A.M.Odlyzko, "*Large deviations of sums of independent random variables*", Acta Arithmetica, **XLIX** (1988), 427-434

A.M.Odlyzko, "*On the distribution of spacings between zeros of the zeta function*", Mathematics of Computation, **48** (1987), 273 – 308

A.M.Odlyzko, www.dtc.umn.edu/~odlyzko/zeta_tables/index.html

M.Rubinstein *zero finder* : http://pmmac03.math.uwaterloo.ca/~mrubinst/L_function_public/L.html

M.Rubinstein and P.Sarnak, "*Chebyshev's Bias*", Experimental Mathematics, **3** (1994), 173-197

T.O. de Silva, "*The first zeros on the critical line of some Dirichlet L series*", www.ieeta.pt

K.Tsang, "*Some $\Omega$ – theorems for the Riemann zeta function*", Acta Arithmetica, **XLVI** (1986), 369-395

38 Southway, Carshalton SM5 4HW, U.K.                    cjmyerscough@btinternet.com